\documentclass[11pt]{amsart}
\def\Max{\operatorname{Max}}
\def\FC{\operatorname{F}}
\def\HF{\operatorname{HF}} 
\def\HP{\operatorname{HP}}
\def\HS{\operatorname{HS}}
\def\deg{\operatorname{deg}}
\def\ker{\operatorname{ker}}

\def\dim{\operatorname{dim}}
\def\depth{\operatorname{depth}}

\def\gin{\operatorname{Gin}}

\newcommand{\m}{\mathfrak m}
\newcommand{\M}{\mathfrak m}
\newcommand{\ZZ}{\mathbb Z}

\newcommand{\QQ}{\mathbb Q}
\newcommand{\R}{\mathcal R}
\newcommand{\B}{\mathcal B}

\def\N{{\bf{N}}}
\def\F{{\mathcal {F}}}

\def\gin{\mathop{\kern0pt\fam0Gin}\nolimits}
\def\gr{\mathop{\kern0pt\fam0gr}\nolimits}
\def\Tor{\mathop{\kern0pt\fam0Tor}\nolimits}
\def\depth{\mathop{\kern0pt\fam0depth}\nolimits}
\def\ini{\mathop{\kern0pt\fam0in}\nolimits}
\def\gcd{\mathop{\kern0pt\fam0gcd}\nolimits}
\def\Grass{\mathop{\kern0pt\fam0Grass}\nolimits}

\def\Rees{{\mathcal R}}
\def\Lex{\operatorname{Lex}}
\def\lex{\operatorname{Lex}}
\def\gin{\operatorname{gin}}

\def\GCD{\operatorname{GCD}}

\newtheorem{theorem}{Theorem}[section]
\newtheorem{lemma}[theorem]{Lemma}
\newtheorem{corollary}[theorem]{Corollary}
\newtheorem{proposition}[theorem]{Proposition}
\newtheorem{remark}[theorem]{Remark}
\newtheorem{example}[theorem]{Example}
\newtheorem{definition}[theorem]{Definition}
\newtheorem{conjecture}[theorem]{Conjecture} 

\begin{document}
\title{Graded rings associated with contracted ideals}
\author{ A. Conca, E. De Negri,\\ A.V. Jayanthan , M. E. Rossi }
\address{A.Conca, E. De Negri,  M. E. Rossi, Dipartimento di Matematica,  Universit\`a di Genova, Via
Dodecaneso 35, 16146 Genova, Italy}
\email{conca@dima.unige.it, denegri@dima.unige.it, rossim@dima.unige.it}
\address{A.V. Jayanthan, School of Mathematics, Tata Institute of Fundamental Research, Homi Bhabha
Road, Colaba, Mumbai, India - 400005.}
\email{jayan@math.tifr.res.in}
\maketitle
\vskip 10mm
\section{Introduction}  The study of the ideals   in a regular local ring $(R, \m)$ of dimension $2$ 
has a long and important tradition dating back to the fundamental work of Zariski \cite{ZS}. More
recent contributions are due to several authors including   Cutkosky, Huneke, Lipman, Sally and Teissier
among others, see \cite{C,C1,H,HS,L,LT}. One of the  main  result in this setting  is  the unique
factorization theorem for complete (i.e., integrally closed) ideals  proved originally by Zariski
\cite[Thm3, Appendix 5]{ZS}. It asserts that  any complete  ideal   can be factorized as a product of
simple complete ideals in a  unique way (up to the order of the factors). By definition, an ideal is 
simple if it cannot be written as a product of two proper ideals.  Another important property of a
complete ideal $I$ is that its reduction number is $1$  which in turns implies that the  associated
graded ring $\gr_I(R)$ is Cohen-Macaulay and  its Hilbert series is well understood;   this is due   
to Lipman and Teissier \cite{LT}, see also \cite{HS}. 

The class of contracted ideals plays an important role in the original work of Zariski as well as in
the work of Huneke \cite{H}.  An
$\m$-primary ideal $I$ of $R$ is contracted if $I=R \cap IR[\m/x]$ for some $x \in \m \setminus \m^2$.
Any complete ideal is contracted but not the other way  round.  The associated graded ring  $\gr_I(R)$ 
of a contracted ideal $I$ need not  be Cohen-Macaulay and its Hilbert series can be very complicated. 

Our goal   is to study  depth,  Hilbert function and  defining equations of the various graded rings
(Rees algebra, associated graded ring and  fiber cone) of homogeneous contracted ideals in the
polynomial ring $R=k[x,y]$ over an  algebraically closed field  $k$  of characteristic  $0$. 

In Section 3 we present several equivalent characterizations of contracted  ideals in the graded and
local case. The main result of this section is Theorem \ref{main-transform}.  It asserts that  the
depth of $\gr_I(R)$ is equal to the minimum of   $\depth  \gr_{I'S_N} (S_N)$, where $S=R[\m/x]$, $I'$
is  the transform of $I$ and $N$ varies in the set of maximal ideals of $S$ containing $I'.$ 

An important invariant of a contracted ideal $I$  of order (i.e., initial  degree) $d$ is the so-called
characteristic  form, which, in the graded setting, is  nothing but $\GCD(I_d)$. Here $I_d$ denotes the
homogeneous component of degree $d$ of $I$.   The more general contracted ideals are those  with  a
square-free  (i.e., no multiple factors) characteristic form.   On the other hand, the more special
contracted ideals    are those  whose characteristic form  is a power of a linear form; these ideals
are exactly the  so-called  lex-segment ideals.   The lex-segment ideals are in bijective
correspondence  with the Hilbert functions (in the graded sense) of graded ideals so that to specify a
lex-segment ideal is equivalent  to specify a Hilbert function.    

In the graded setting Zariski's factorization theorem for  contracted ideals  \cite[Thm1, Appendix
5]{ZS} says that any contracted ideal $I$ can be written as a product of lex-segment ideals
$I=\m^cL_1\cdots L_k$.  Here each $L_i$ is a lex-segment monomial ideal with respect to an  appropriate
system  of coordinates  $x_i, y_i$ which depends on
$i.$ Furthermore  $L_i$  has exactly one generator in its initial degree which is a power of $x_i$. 

As a consequence of Theorem \ref{main-transform} we have  that the depth of $\gr_I(R)$ is equal to the
minimum of the  depth of
$\gr_{L_i}(R) $ (see Corollary \ref{L1}). We can also express the Hilbert series of $I$ in terms of the
Hilbert series of the $L_i$'s and of the characteristic form of $I$ (see Proposition \ref{hfperI}).  

For a contracted ideal with a square-free characteristic form we show in Theorem \ref{gcdsqfree} that
the Rees algebra $\Rees(I)$, the associated graded ring $\gr_I(R)$ and the fiber cone $\FC(I)$ are all
Cohen-Macaulay with expected defining equations in the sense of
\cite{Vas} and \cite{MU}.  Furthermore  $\Rees(I)$ is normal, the fiber cone $\FC(I)$ is reduced and
we  determine  explicitly  the Hilbert function  of $\gr_I(R)$. 

Section 3 ends with a statement and a conjecture on the coefficients of the $h$-vector of a contracted
ideal. Denote by $h_i(I)$ the $i$-th coefficient  of the $h$-vector of $I$ and by $ \mu(I)$ the minimal
number of generators of $I.$  We show that for any contracted (or monomial)  ideal $I$  one has
$h_1(I)\geq (\mu(I)+1)\mu(I)/2$  (Proposition \ref{boundh1}) and we conjecture that
$h_2(I)\geq 0$.

Sections 4 and 5 are devoted to the study of the lex-segment ideals. This class is important since, as
we said above, information about the associated graded ring of a  contracted ideal $I$  can be derived
from information about the associated graded rings of the lex-segment ideals appearing in the 
Zariski's factorization of $I.$  Another reason for studying the associated graded rings of lex-segment
ideals comes from Section 2.  There it is proved that if $I$ is any ideal and
$\ini(I)$ is its initial ideal with respect to some term order then
$H_I(n)\leq H_{\ini (I)}(n)$ for all $n$ provided $\depth \gr_{\ini (I)}(R)>0$. In two  variables 
initial ideals in generic coordinates are lex-segment ideals. We  detect several classes of lex-segment
ideals for which the associated graded ring is Cohen-Macaulay or at least has positive depth. In
particular, consider the lex-segment ideal $L$  associated with  the Hilbert function of an ideal
generated by generic forms $f_1,\dots, f_s$; equivalently, set $L=\ini((f_1,\dots, f_s))$, where the
forms $f_i$ are generic  and in generic coordinates. We show that  $\depth\gr_{L}(R)>0$ (see Theorem
\ref{positivedepth1}). Furthermore  $\gr_{L}(R)$ is Cohen-Macaulay if  the forms $f_i$ have all the
same degree.  In Section 6 we describe the  defining equations of the Rees algebra of various classes
of lex-segment ideals. 

Some of the results and the examples  presented in  this  paper  have been inspired and suggested by
computations   performed  by   the  computer algebra system CoCoA   \cite{Co}. In particular we have
made extensive use of the local algebra package.  

\vskip 5mm The authors thank G.Valla for many useful discussions on the topic of this paper. Part of
this work was done while the third author was visiting University of Genova.  He  would like to thank 
the Department of Mathematics for the support and the  nice hospitality.

\section{Initial ideals  and  associated graded rings} Let $R $ be a regular local ring of dimension
$d$, with maximal ideal
$\M$ and residue field $k$, or, alternatively, let
$R=k[x_1,\dots,x_d]$ be a polynomial ring over a field $k$, and $\M$ a maximal ideal  of $R.$
Throughout the paper we  assume that $k=R/\M $  is algebraically closed  of characteristic 0. 
Moreover, let $I$ be an
$\M$-primary ideal.  For every integer $n, $ the length $\lambda ( R/I^{n+1}  ) $ of $ R/I^{n+1} $ as
$R$-module is finite. For  $n\gg 0,
$  $\lambda ( R/I^{n+1}  ) $  is a polynomial $\HP_I(n)$ of degree $d$ in $ n $. The polynomial
$\HP_I(n)$ is called the Hilbert-Samuel polynomial of $I$ and one has $$ \HP_I(n)= \frac {e(I)}{d!} n^d
+ \ \ lower\ \ degree \
\ terms. \ $$ In particular, $e(I) $ is the ordinary multiplicity of the associated graded ring
$\gr_I(R)$ to $I$,
$$
\gr_I(R)= \oplus_{n\ge 0} I^n/I^{n+1}. 
$$

The Hilbert function $\HF_I(n)$ of $I$ is defined as
$$
\HF_I(n)= \lambda(I^n/I^{n+1} )
$$ and it is by definition the Hilbert function of $\gr_I(R).$ The Hilbert series of $I$ is 
$$
\HS_I(z)= \sum_{n\ge0} \HF_I(n)z^n = \frac {h_0(I)+h_1(I)z+\ldots+h_s(I)z^s}{(1-z)^d},
$$ where $h_0(I)=\lambda(R/I)$ and $\sum_{i=0}^s h_i(I)=e(I).$

\vskip 5mm In the local case, most important tools for studying the associated graded ring are  minimal
reductions  and superficial elements.  Those tools are not available in the non-local case, so we need
to pass to the localization.
\begin{lemma} \label{estensioni} {  Let $S$ be  a flat extension of a ring $R$ and let $I\subset R$ be
an ideal. Suppose that $S/IS\simeq R/I$ as $R$-modules. Then
$$\gr_I(R)\simeq \gr_{IS}(S).$$ }
\end{lemma}

\begin{proof} It is enough to prove that $I^n/I^{n+1}\simeq I^nS/I^{n+1}S$ as
$R$-modules. Since $S$ is a flat extension of $R$, one has
\begin{eqnarray*}
  I^nS/I^{n+1}S & \simeq & I^nS\otimes_S S/IS\simeq I^n\otimes_R
  S\otimes_S S/IS \\
  & \simeq & I^n\otimes_R S/IS\simeq I^n\otimes_R
  R/I\simeq I^n/I^{n+1},
\end{eqnarray*} and the multiplicative structure is the same.
\end{proof}
\bigskip

\begin{remark} \label{casi} {\rm{ We  apply the above lemma to our setting.  The localization
$R_\M $ is a flat extension of $R  $ and $R/I \simeq R_\M/I_\M.$  By the above lemma we get
$$\gr_I(R) \simeq \gr_{I_\M}(R_\M) $$

In particular, one has $\HF_{I_\M}(n) =\HF_{I }(n )  $  and hence
$\HS_{I_\M}(z) =\HS_{I }(z ).  $}}
\end{remark}
\vskip 5mm

We are interested in studying the behaviour of  $\HF_I(n)$ under Gr\"obner deformation.  In the
following we denote by $I$ an
$\M$-primary ideal of $R=k[x_1,\dots,x_d]$ with $\M=(x_1,\dots,x_d) $. We fix a term order on $R$ and
consider  the initial ideal $ \ini(I) $ of $I$.  Recall that  
$\lambda(R/I) =\lambda( R/\ini(I)).$ 

We want to  compare $\HF_I(n) $ and $\HF_{\ini(I)}(n). $  First note that $$ e(I) \le e(\ini(I)). $$
This inequality  follows easily from the fact that the multiplicities can be read from the leading
coefficients of the Hilbert-Samuel polynomials of $I$ and of $\ini(I)$.  In fact, since
$\ini(I)^n \subseteq \ini(I^n)$, we have
$$
\lambda (R/I^n) =\lambda(R/\ini(I^n)  )\le\lambda(R/\ini(I)^n)\eqno(1)
$$ and hence for $n \gg 0$ we get the required inequality on the multiplicities. Notice that   in
\cite[4.3 ]{DTVVW}   the equality has  been characterized .

As a consequence of the next lemma one has  that the same inequality holds for the Hilbert function
asymptotically, in a more general setting. Moreover, under some assumption, the inequality holds from
the beginning.

\begin{lemma} \label{lemmadepth}  {  Let $J$ be an $\m$-primary ideal in $R=k[x_1,\dots,x_d]$ and let
$\F=\{\F_n\}_{n\ge 0}$ be a filtration of ideals, such that
$J^n\subseteq \F_n$ for every $n\ge 0$. Then
\begin{itemize}
  \item[(1)] $\lambda(\F_n/\F_{n+1})\le \lambda(J^n/J^{n+1})$ for
    $n\gg 0$;

  \item[(2)] If $ \depth  \gr_{J}(R)>0$, then
    $\lambda(\F_n/\F_{n+1})\le \lambda(J^n/J^{n+1})$  for every $n\ge
    0$.
\end{itemize} }
\end{lemma}

\begin{proof}  (1) Since for every $n\ge 0$ we have  $J^n\subseteq \F_n$, it is equivalent to prove
$$ \lambda(\F_n/J^{n}) \le \lambda(\F_{n+1}/J^{n+1}).$$
\noindent By Remark \ref{casi}, and since
$\lambda(\F_n/\F_{n+1})=\lambda(\F_nR_{\M}/\F_{n+1}R_{\M})$, we may transfer the problem to the  local
ring $R_{\M}$, identifying $J$ with $JR_{\M}$ and $\F_n$ with $\F_nR_{\M}$.  Let $a $ be a superficial 
element for $J$ and consider the following exact sequence induced by the multiplication by $a$:
$$ 0 \to [(J^{n+1}  : a ) \cap \F_n  ]/J^{n}  \to \F_n/J^n \stackrel{\cdot a}{\to}
\F_{n+1}/J^{n+1} \to \F_{n+1}/a\F_n+J^{n+1} \to 0.$$ Since $a$ is superficial and regular, one has
$J^{n+1}  : a=J^n$ for
$n\gg 0$, and this proves (1).

(2) If depth $ \gr_{J}(R) >0, $ then $\overline a \in J /J^2  $ is regular (see \cite[2.1]{HM1}) and $
\ \ J^{n+1}  :~ a~ = J^n  $ for every
$n. $ This forces the required inequality and concludes the proof.
\end{proof} As a consequence of the above lemma one has:

\begin{theorem} \label{depth} Fix any term order on $R=k[x_1,\dots,x_d]$, and let $I$ be an
$\m$-primary ideal in $R$. The following facts hold:
\begin{itemize}
  \item[(1)] $\HF_I(n)\le \HF_{\ini(I)}(n)$ for $n\gg 0$;

  \item[(2)] If $\depth\gr_{\ini(I)}(R)>0$, then $\HF_I(n)\le
    \HF_{\ini(I)}(n)$  for every $n\ge 0$.
\end{itemize}
\end{theorem}

\begin{proof} Since $\lambda(R/I^n)= \lambda(R/\ini(I^n)) $ for every $n,$ one has
$$\HF_I(n) = \lambda(I^n/I^{n+1})=\lambda(\ini(I^n)/\ini(I^{n+1})).$$ Now the results follow by
applying Lemma \ref{lemmadepth} with
$J=\ini(I)$ and $\F_n=\ini(I^n)$. Note that part (1)  can also be proved directly from equation $(1).$ 
\end{proof}
\bigskip A lex-segment ideal is a monomial ideal $L$  such that  whenever $n,m$ are monomials of the
same degree with  $n>m$ in the lexicographical order then $m \in L$ implies $n\in L$. Macaulay's
Theorem on Hilbert function implies  that for every homogeneous ideal $I$ there is a unique lex-segment
ideal $L$ with $\dim I_s=\dim L_s$ for all $s$. We will denote this ideal by $\lex(I)$. Note however
that $\lex(I)$ depends only on the Hilbert function of $I$.  

The following examples show that the conclusion of part (2) in Theorem \ref{depth} does not hold if the
depth of
$\gr_{\ini(I)}(R) $ is 0.
\vskip 5mm
\begin{example}\label{ex11} {\rm {(a) Consider $ R=\QQ[ x,y ]$ equipped with the lexicographic order,
with $x>y$.  If $ I= (x^5, x^4y^2, x^2y^5(x+y), xy^8, y^{10})$, then $\ini(I)= (x^5, x^4y^2,
x^3y^5,x^2y^7, xy^8, y^{10}).$ In this case the associated graded ring to $\ini(I)$ has depth zero and
one has:
$$
\HS_I(z)=\frac{32+14z+6z^2-2z^3}{(1-z)^2}, \qquad
\HS_{\ini(I)}(z)=\frac{32+16z+4z^2-2z^3}{(1-z)^2}.
$$ Thus  $ \HF_I(2)=130>128=
\HF_{\ini(I)}(2) $ in spite of the equality $e(I)=e(\ini(I))$, and of the fact that
$\HF_I(n)=\HF_{\ini(I)}(n)$ for $n\ge 3$.  Note that
$\ini(I)$ is a lex-segment ideal, thus in particular, it is also the generic initial ideal of $I$.

\medskip

(b) Let $I=(x^9, x^7y, x^6y^3,  x^5y^5, x^4y^{12},  x^3y^{13}, x^2y^{14},  xy^{17}, y^{19})\subseteq
\QQ[x,y]$. Its generic initial ideal is the lex-segment ideal $L=(x^8, x^7y^2, x^6y^3, x^5y^5,
x^4y^{12}, x^3y^{13}, x^2y^{14}, xy^{17}, y^{19})$.  Notice that $I$ is contracted, since it has the
same number of generators of $L$ (contracted ideals will be defined and studied in Section
\ref{contracted}$)$.  Neither  in this case the inequality holds. In fact one has
$$
\HS_I(z)=\frac{85+42z+10z^2-3z^3}{(1-z)^2} \qquad \mbox{and} \qquad
\HS_L(z)=\frac{85+43z+7z^2-z^3}{(1-z)^2},
$$ thus  $ \HF_I(2)=349>348= \HF_{L}(2) $. Also in this case
$\depth\gr_{\ini(I)}(R)=0$.}}
\end{example}
\vskip 5mm

By Theorem \ref{depth}, in $k[x,y]$, one has $$
\lambda(I^n/I^{n+1})\leq \lambda(\Lex(I)^n/\Lex(I)^{n+1}) \mbox{\ \ \
\  for every }n\gg 0.$$ Such  inequality  does not hold in $3$ or more variables, see the next example.
\vskip 3mm

\begin{example} \label{3var} {\rm{ Let $I=(x^2,y^2,xy, xz^2, yz^2 , z^4)\subset \QQ[x,y, z ]$.  One
has\\ $\Lex(I)=(xz, xy, x^2, yz^2, y^2z, y^3, z^4),$ and
$$\HS_I(z)=\frac{8+8z}{(1-z)^3} ,\ \ \ \ \HS_{\Lex(I)}(z)=
\frac{8+7z}{(1-z)^3}.$$ Thus $ \HF_I(n) \ge \HF_{\Lex(I)}(n) $ for every $n\ge 1$, and also $e(I) >
e(\Lex(I)). $ }}
\end{example}
\bigskip

In the last part of the section we restrict ourselves to the case of dimension two, and we collect some
Cohen-Macaulayness criteria for the associated graded ring. We recall  the following important
result.

\begin{proposition} \label{dim2}\cite[Thm. A, Prop. 2.6]{HM}. Let $I$ be an $\M$-primary ideal of a 
regular local ring $(R,\m)$ of dimension two and  $J$  a minimal reduction of $ I.$ Then  $\gr_I(R)$ is
Cohen-Macaulay if and only if $I^2=JI.$
\end{proposition}

It follows from \cite[5.5]{LT} and  \cite[3.1]{HS} that:

\begin{proposition}
\label{ic} Let $I$ be an  $\M$-primary ideal in a regular local ring $(R,\m)$ of dimension two. If $I$ 
is  integrally closed, then $ \gr_I (R)$ is Cohen-Macaulay.
\end{proposition}

Proposition \ref{dim2} does not have a corresponding  version in the graded case since minimal reductions
need not to exist in that setting. But the next corollary holds both in the graded and in the local
setting.

\begin{proposition} \label{HScorta} Let $I$ be an  $\M$-primary ideal in a regular ring $ R $ of
dimension two. Then  $\gr_I(R) $ is Cohen-Macaulay if and only if
$$
\HS_I(z)=\frac{h_0(I)+h_1(I)z}{(1-z)^2}.
$$
\end{proposition}
\begin{proof} This is a simple consequence of the fact that $\gr_I(R) \simeq \gr_{I_{\M}}(R_{\M}). $  
In fact, by Proposition \ref{dim2}, $\gr_{I_{\M}}(R_{\M}) $ is Cohen-Macaulay if and only if $I_{\M}^2
=J I_{\M}  $ for some reduction $J$ of
$I_{\M}.$ Since $R_{\M}$ is a local Cohen-Macaulay ring, this is equivalent to
$$ \HS_{{I_{\M}}}(z)= \frac{h_0(I)+h_1(I) z}{(1-z)^2}$$  (see for example \cite[2.5]{GR}). The
conclusion follows by Remark \ref{casi}.
\end{proof}

\vskip 3mm

 We will apply the above criteria for proving the Cohen-Macaulayness of the associated graded rings of
certain classes of monomial ideals.

\vskip 3mm

Let  $R=k[x,y]$  and denote by $\m$ the ideal $(x,y)$. There are many ways of encoding an $\m$-primary
monomial ideal $I$.  Among them we choose the following.

Set  $d=\min\{ j : x^j \in I\}$. Then for $i=0,\dots, d$  we set
$a_i(I)=\min\{ j : x^{d-i}y^j \in I\}$ and $$a(I)=(a_0(I), a_1(I),
\dots, a_d(I)).$$ The sequence $a(I)$ is said to be the column sequence of $I$. By the very definition
we have that $a_0(I)=0$ and
$1\leq a_1(I)\leq a_2(I)\leq \dots \leq a_d(I)$. Conversely, any sequence satisfying these conditions
corresponds to a monomial ideal. For example
\begin{eqnarray*}
  I = (x^3, xy^3, y^5) & \longleftrightarrow & a(I) = (0, 3, 3, 5), \\
  I = (x^4, x^3y, x^2y^4, xy^7, y^9) & \longleftrightarrow & a(I) =
  (0, 1, 4, 7, 9).
\end{eqnarray*} It is easy to see that: 
$$\lambda(R/I)=|a(I)|=\sum_i a_i(I).$$

Note also that the minimal generators of $I$ are the monomials
$x^{d-i}y^{a_i(I)}$  with $a_i(I)<a_{i+1}(I)$ or $i=d$. 

\begin{remark} \label{lexiscon} {\rm In two variables, as we are, the  $\m$-primary lex-segment ideals
correspond  exactly to strictly increasing  column sequences. In other words any lex-segment ideal $L$
can be written  as $$L=(x^d, x^{d-1}y^{a_1},\dots,   xy^{a_{d-1}}, y^{a_d}),$$ where
$0=a_0<a_1<\dots < a_d$. In particular, the minimal number of generators of $L$ is exactly one more
than the initial degree. Ideals with this property are called contracted, see Section
\ref{contracted}.  Furthermore, in characteristic $0$, the lex-segment ideals are exactly the Borel
fixed ideals and this implies that the generic initial ideal 
$\gin(I)$ of $I$ is equal to $\lex(I)$.  }
\end{remark}  

\medskip

The $b$-sequence (or differences sequence) of $I$ is  denoted by
$b_1(I),\dots,b_d(I)$ and defined as
$$b_i(I)=a_{i}(I)-a_{i-1}(I).$$

We will use $a_i$ for $a_i(I)$ and $b_i$ for $b_i(I)$ if there is no confusion.

If $I$ and $J$ are monomial $\m$-primary ideals,  then the column sequence of the product $IJ$ is 
given by 
$$a_i(IJ)= \min \{ a_j(I)+a_k(J)  : j+k=i\}.$$ In particular, the column sequence of $I^n$ is given by: 
$$a_i(I^n)=\min\{ a_{j_1}(I)+\cdots+a_{j_n}(I) : j_1+\cdots+j_n=i\}.$$
\vskip 3mm
\begin{example}\label{differenze} {\rm{ Let $I$ be a monomial ideal with $b$-sequence $b_1,\dots, b_d$. 
\begin{itemize}
  \item[(a)] Assume  $b_1\leq b_2\leq \dots \leq  b_d$.  Then it is
    easy to see that for every $n\in\N$ one has  $a_i(I^n)=(n-r)a_q(I)+ra_{q+1}(I)$, where
    $i=qn+r$ with $0\leq r<n$.  Summing up, we have
    $|a(I^n)|=n^2|a(I)|-{n \choose 2} a_d(I)$. It follows that 
    $$\HS_I(z)=\frac{\lambda(R/I)+(\lambda(R/I)-a_d(I))z}{(1-z)^2}.$$
    Note that  by \cite[Ex. 4.22]{E} the ideal $I$ is integrally
    closed.
    \vskip 2mm
  \item[(b)] Assume  $b_1\geq b_2\geq \dots\geq  b_d$.  Then it is
    easy to see that  for every $n\in\N$ one has $a_i(I^n)= q a_d(I)+ a_{r}(I)$, where $i=qd+r$
    with $0\leq r<d$.  Summing up, we have $|a(I^n)|=n|a(I)|+{n
    \choose 2} d a_d(I)$. It follows that
    $$\HS_I(z)=\frac{\lambda(R/I)+(da_d(I)-\lambda(R/I))z}{(1-z)^2}.$$
    Moreover, the ideal $J=(x^d,y^{a_d})$ is a minimal reduction of $I$, and it holds $I^2=JI$.
  
\end{itemize} }
  In both cases the associated graded ring to $I$ is Cohen-Macaulay
    by \ref{HScorta}. }

\end{example} 
\medskip

\section{Contracted ideals }
\label{contracted}

Let $R$ be either a polynomial ring   over a field or a regular local ring. Denote by $\m$ the
(homogeneous) maximal ideal of $R$. Assume the residue field $k$ is algebraically closed of
characteristic 0, and the dimension of $R$ is $2$.  Let $a \neq 0 $ be an element of
$\M, $ the order $o(a)$ of $a$ is the $\M$-adic valuation  of $a,$ that is, the greatest integer $n$
such that  $a\in \M^n$. If $o(a)=r,$ then we denote by $a^* $ the initial form of $a$ in $\gr_{\M}(A),$
that  is $a^*=\overline a \in \M^r/\M^{r+1}.$

Let $I$ be an $\M$-primary ideal, homogeneous in the graded case. Denote by $\mu(I)$ the minimum 
number of generators of $I$ and by
$o(I)$ the order of $I$, that is, the  largest $h$ such that
$I\subseteq \m^h$. In the graded case $o(I)$ is simply the least degree of  non-zero elements in $I$.
In the local case, if $I^* $ is the homogeneous ideal of $\gr_{\M}(R) $ generated by the initial forms
of the elements of $I$, then $o(I)$ is the least degree of an element in $I^*$.  As in \cite{ZS}, we
call {\it characteristic form} the GCD of the elements of degree $o(I)$ in $I^*.$ In the graded case the
characteristic form is just the GCD of the elements of degree $o(I)$ in $I.$

By the Hilbert-Burch theorem, $I$ is generated  by the maximal minors of a $(t-1) \times t$ matrix, say
$X$, where $t=\mu(I)$. It follows that $\mu(I)\le o(I)+1$.

\begin{remark} \label{degmatta} {\rm { In the graded  setting, if  $g_1\leq \dots \leq g_{t}$ are the
degrees of the generators of $I$ and $s_1\leq \dots \leq s_{t-1}$ the degrees of the syzygies, then the
$ij$-entry of $X$ has degree $s_i-g_j$. Here we use the convention that $0$ has any degree. The matrix
$(u_{ij})$,
$u_{ij}=s_i-g_j$,   is called the degree matrix of $I$. It is easy to see that $u_{ij}$ must be
positive for all $i,j$ with $j-i\leq 1$ and that $o(I)=\sum_{i=1}^{t-1}  u_{i,i+1}$.}}
\end{remark} 

We have the following: 
\begin{proposition} 
\label{contrattiscaduti} Let $I$ be  an $\M$-primary ideal.  The following conditions are equivalent: 
\begin{itemize} 
  \item[(1)] $\mu(I)=o(I)+1$,

  \item[(2)] there exists   $\ell\in \m \setminus  \m^2$ such that
    $I:(\ell)=I:\m$,

  \item[(3)]  there exists   $\ell\in \m \setminus  \m^2$ such that $I=IS\cap
    R,$ where $S=R[\m/\ell]$. 
\end{itemize} 

Furthermore, in the graded case,   $(1),(2), (3)$ are equivalent to any of the following conditions: 

\begin{itemize} 
  \item[(4)] $I$ is Gotzmann, i.e., $I$ and $\Lex(I)$ have the same number of
    generators, 

  \item[(5)] $I$ is componentwise linear, i.e., the ideal generated by every
    homogeneous component of $I$ has a linear resolution, 

  \item[(6)] for every $h\geq o(I)$, the degree $h$ component of $I$
    has the  form $I_h=f_h R_{h-s_h},$ where $f_h$ is a homogeneous
    polynomial  of degree $s_h$,

  \item[(7)]  $u_{i,i+1}=1$ for $i=1\dots, \mu(I)-1$.  
\end{itemize} 
\end{proposition} 

\begin{proof}  The equivalence between (1),(2) and (3) is proved in \cite[2.1,2.3]{H} in the local case
and the arguments work also in the graded case.  In the graded case the equivalence between (1) and (7)
follows from Remark
\ref{degmatta}.  The equivalence of  (1)  and (4) holds because, obviously, $o(I)=o(\Lex(I))$ and any
lex-segment ideal in $2$ variables satisfies (1), see Remark \ref{lexiscon}. That  (4) implies (5)  is
a general fact  \cite[Ex.1.1b]{HH} while that (5) is equivalent to  (6) follows from the fact that in
two variables the only ideals with linear resolution have the form $f\m^u,$ where $f$ is a form. To
conclude, it suffices to show that (6) implies (2), where $\ell$ is any linear form not dividing the
characteristic form of $I$ and this is  an easy check. 
\end{proof}
\begin{definition} 
\label{defcontratti} An $\m$-primary ideal $I\subset R$  is said to be contracted if it satisfies the
equivalent conditions of Proposition \ref{contrattiscaduti}. In particular, we say that $I$ is
contracted from  $S=R[\m/\ell]$ if
$I=IS\cap R$.
\end{definition}

If $I$ is contracted, then conditions (2) and  (3) hold true for any
$\ell$ such that its initial form in $\gr_\m(R)$ does not divide the characteristic form of $I.$  In
the graded case, this just means that
$\ell$ is a linear form not dividing the $\GCD$ of the elements of degree $o(I)$ of $I$. Such an
element $\ell$ is called  {\it {coprime}} for $I$. Since $k$ is infinite,  coprime elements   for $I$
exist. More generally given a finite number of ideals one can always find an element which is coprime
for any  ideal.

By Remark \ref{lexiscon}, the minimum number of generators of any lex-segment ideal $L$  is exactly one
more than the initial degree; hence $L$ is  contracted. Moreover,  $ y+ ax$ is coprime for $L$ for all
$a \in k.$

\medskip
\begin{remark} \label{lambda} {\rm{The homogeneous component $I_h$ of  a  homogeneous contracted ideal
$I$  has the form   $I_h=f_h R_{h-s_h}$ for all
$h\geq o(I)$.  The element $f_h$ is the $\GCD$ of the elements  in
$I_h$. Furthermore it divides $f_{h-1}$ for all $h>o(I)$.  Here
$s_h=\deg f_h$ is also the dimension of $R/I$ in degree $h$.  So we have 
$$\lambda(R/I)={ \mu(I) \choose 2} +\sum_{h\geq o(I)}  s_h.$$ The number of generators of $I$ in degree
$h> o(I)$ is
$s_{h-1}-s_{h}$. The lex-segment ideal $L=\Lex(I)$ associated with $I$ has the following form
$L_h=x^{s_h} R_{h-s_h}$. }}
\end{remark} 

Next we give a characterization  of contracted ideals in terms of the Hilbert-Burch matrix: 

\begin{proposition}
\label{newcontracted}  Let $R=k[x,y]$ and let $d$ be a positive integer. Let
$\alpha_1,\dots, \alpha_d$ be elements of the base field $k$ and
$b_1,\dots, b_d$ positive integers. Then the ideal generated by the
$d$-minors of the $d\times (d+1)$ matrix 
$$
\left (
\begin{array}{cccccccccccccccccccccc} 
  y^{b_1}   &  x+\alpha_1y \ &    0        & 0       & \cdots & \cdots & 0 
  \\
  0     &     y^{b_2}         &x+\alpha_2y\ & 0       & \cdots & \cdots &
  0   \\
  \cdots         &\cdots       &  \cdots   & \cdots  & \cdots & \cdots&
  \cdots\\
  \cdots         &\cdots       &  \cdots   & \cdots  & \cdots & \cdots&
  \cdots\\
  0              &\cdots            &\cdots          &\cdots &0  & 
  y^{b_d}  & x+\alpha_dy  
\end{array}
\right)
$$ is contracted of order $d$, and  $y$ is coprime for $I$. Conversely,
   every contracted ideal $I$ of order $d$ can be realized, after  a change of coordinates, in this
way. 
\end{proposition} 

\begin{proof} That the ideal of minors of such a matrix is contracted follows directly from the
definition. Conversely, assume that $I$ is contracted and let $L$ be its associated lex-segment ideal.
Say
$a=(a_0,\dots, a_d)$ is the column sequence of $L$.    The matrix above with   $b_i=a_i-a_{i-1}>0$ and 
all the $\alpha_i=0$  defines
$L$. On the other hand, for every choice of the $\alpha_i$ we get a contracted ideal with the Hilbert
function of $I$. We will show that a particular choice of the $\alpha_i$ will define the ideal $I$. By
definition, $I$ is determined by $d$ and by the form $f_h=\GCD (I_h)$ for $h\geq d$.  Since we assume
that $k$ is algebraically closed, every $f_h$ is a product of linear forms. Since $f_{h+1}$ divides $
f_{h}$, we may find linear forms $\ell_{d-s_h+1} ,\dots, \ell_d$ so that $f_h=\Pi_{j=d+1-s_h}^d
\ell_j,$ where $s_h=\deg f_h$. Take linear forms  $\ell_1,\dots, \ell_{d-s_h}$ in any way. Then take a
system of coordinates $x,y$ so that  $y$ does not divide any of the $\ell_i$. In other words, up to
irrelevant  scalars, $\ell_i=x+\alpha_iy$ for all
$i=1,\dots, d$. We claim that this choice of the $\alpha_i$ works. The degree of $f_h$ is determined by
the Hilbert function.  So,  it is enough to show that $f_h$ divides every $d$-minor of degree $\leq h$
of the matrix. This is easy to check.
\end{proof} 

\vskip 3mm

An important property of contracted ideals is the following: 

\begin{proposition} 
\label{powconiscon} The product of contracted ideals is contracted. 
\end{proposition} 

\begin{proof} The proof  given in the local case in \cite[2.6]{H} works also in the graded setting.
\end{proof}

\medskip Next we recall Zariski's factorization  theorem for contracted ideals
\cite[Thm1, Appendix 5]{ZS}: 

\begin{theorem} \label{Zariski2} Let $I$ be a contracted ideal of order $d$ and characteristic form $g$
of degree $s$.  Let $g=g_1^{\beta_1}\dots g_k^{\beta_k}$ be the factorization of $g,$ where the
$g_i $'s  are distinct  irreducible forms. Then $I$ has a unique factorization as:
$$I=\m^{d-s}L_1L_2\cdots L_k,$$ where the $L_i $'s are contracted ideals with characteristic form
$g_i^{\beta_i}$. 
\end{theorem}

Since we assume that $k$ is algebraically closed the $g_i $ 's are indeed linear forms. In the graded
setting it follows that the $L_i $ 's are lex-segment ideals in a system of coordinates with $g_i$ as
first coordinate.  The ideals $L_i$'s  can be described quite explicitly in terms of the data of $I$: if
$g_i$ appears in the $\GCD$ of $I_{d+j}$ with exponent
$\beta$, then the $\GCD$ of $L_i$ in degree $\beta_i+j$ is exactly
$g_i^\beta$. Then Remark \ref{lambda} implies: 

\begin{lemma} 
\label{lunperI} With the notation of Theorem \ref{Zariski2} we have: 
$$\lambda(R/I)=\sum_{j=1}^k \lambda(R/L_j)+{ d+1 \choose 2}-\sum_{j=1}^k { \beta_j+1 \choose 2} . $$
\end{lemma} 

From the factorization of $I$  given in  Theorem \ref{Zariski2} immediately follows that 

$$I^n=\m^{n(d-s)}L_1^nL_2^n\cdots L_k^n$$ is the analogous factorization for $I^n$. Applying Lemma
\ref{lunperI} to $I^n$, summing up and using the formula 

$$ \sum_{n=0}^{\infty} { (n+1)\gamma +1 \choose 2} z^n=\frac{\gamma(1-z)+\gamma^2(1+z)}{2(1-z)^3}$$  
we obtain: 

\begin{proposition} 
\label{hfperI} With the notation of Theorem \ref{Zariski2} we have: 
$$HS_I(z)=\sum_{j=1}^k HS_{L_j}(z)+\frac{{d+1 \choose 2}+{d\choose 2}z-  \sum_{j=1}^k \left[{\beta_j+1
\choose 2}+{\beta_j \choose 2}z\right]}{(1-z)^2 }
$$ and in particular,
$$e(I)=\sum_{j=1}^k e(L_j)+d^2-\sum_{j=1}^k \beta_j^2.$$  
\end{proposition}

Similarly one can write all the coefficients of the Hilbert-Samuel polynomial of $I$ in terms of those
of the  $L_i$'s.

\bigskip

In this part of the section $(R,\m)$ will denote a regular local ring of dimension two and  $I$ an
$\m$-primary ideal.  We consider a coprime element  $\ell$ for $I,$ and  we fix a minimal system of
generators of
$\M=( x,\ell ) $.
\vskip 5mm
\noindent  We define now the {\it transform} of an ideal $I$ (not necessarily contracted)  in
$S=R[\M/\ell]$.  If $a$ is in $I$ and $d=o(I)$ is the order of $I$, then $a/\ell^d$ is in $S$  and we 
may write
$$IS=\ell^dI',$$  where $I'$ is an ideal of $S$.  Such an ideal $I'$ is called the transform of $I$ in
$S$.

Notice  that if $I=(f_1,\dots,f_t)$, then
$I'$ is generated (not minimally in general) by $f_1/\ell^d,\dots,f_t/\ell^d$. 
 If $I$ and $J$ are two ideals of $R$, then $(IJ)'=I'J'$. In fact, if $d=o(I)$ and $s=o(J)$, then
$(IJ)S=(IS)(JS)=\ell^d \ell^s I'J'$. Therefore
$(IJ)S=\ell^{d+s}(IJ)'$, so the conclusion follows. In particular,
$(I^n)'=(I')^n$ for any integer $n.$

If $I=\M^d, $ then $I'=S.$ In particular, if $I=\M^sJ,$  then  $I'=J'
$  in $S. $ 

\vskip 5mm In the following   we   always denote by $ I'$ the transform of $I$ in
$S=R[\M/\ell]$.
\vskip 5mm

We note that $S= R[\M/\ell]=R[x/\ell]$ is  isomorphic to  the ring
$R[z]/(x-z\ell) $. The ring $S$ is not  local and its maximal ideals
$N$   which contain $\M $ are in one to one  correspondence with the irreducible polynomials $g$ in
$k[z].$

We denote by $T$  the localization of $S$ at one of its maximal ideal
$N$. Then $T$ is a $2$-dimensional regular local ring called the {\it first quadratic transform } of
$R.$ In algebraic geometry this construction is the well-known ``locally quadratic" transformation of
an algebraic surface, with center at a given simple point $P$ of the surface.

If $I'$ is the transform  of $I$ in $S$, then $(I')_N=I'T $ is the transform of $I$ in $T $ and
$IT=\ell^d (I')_N $ if $d=o(I).$ We remark   that $\ell$ is a regular element both in $S$ and $T$.  It
is known that if $I $ is primary for $\M $ and $N$ is a maximal ideal which contains $I',$ then $(I')_N
$ is primary for $N $  or is a unit ideal. If $I$ is contracted, then $(I')_N $ is a unit ideal if and
only if $I=\M^d $ (see \cite[Prop. 2 and Corollary, Appendix 5]{ZS}). \vskip 2mm \par \noindent In the
following we assume $I$ is not a power of the maximal ideal.  \vskip 2mm \par Then $I'$ is a
zero-dimensional ideal of $S$ which is not necessarily primary. We will denote by Max($I'$) the set of
the maximal ideals associated to
$I'.$ The maximal ideals in  Max($I'$)  depend on the characteristic form of $I $ and on the field $k.$

Denote by $T$ any localization of $S $ at a maximal ideal $N \in $ Max($I'$).  The following easy facts
will be useful in the proof of Theorem \ref{main-transform}.
\vskip 3mm
\begin{remark} \label{reduction}  {\rm{ Let $d=o(I), $  $\ell$ be  coprime for $I$, $J$ be a minimal
reduction of $I$ and  $S=R[\M/\ell] $. The following facts hold:
\begin{itemize}
  \item[(1)] $\ell$  is  coprime for $I^n $ and for $\M^sI^n$ for
    every positive integers $n$ and $s.$ In particular, if $I$ is
    contracted from $S$, then  $I^n$ and  $\M^sI^n$ are contracted from
    $S$.
 
  \item[(2)] $\ell$ is  coprime for $J.$  In fact, there
    exists $n$ such that $I^{n+1}=JI^n $ and since  $ o(J)=d $, to conclude  it is
    enough to look at the minimal degree $nd+d $ part of the
    corresponding  ideals of the initial forms.

  \item[(3)] If $I $ is contracted, then $JI $ is  contracted
    (\cite{H}, proof of Thm. 5.1), and  by  {\rm (2)} it is
    contracted from  $S$.
 
  \item[(4)]If $J=(a,b)  $  ,  then
    $J'=( a/\ell^d, b/\ell^d)  $ is a minimal reduction of $I'$ both
    in $S $ and $T.$   In fact, if $I^{n+1}=JI^n,$ then $I^{n+1}S=JI^n
    S.$ Since $o(I^{n+1})=o(JI^n)$, we have  $$ (I')^{n+1}=J'
    (I')^n.$$ In particular, from the last equality it follows easily
    that Max($I'$)  =   Max($J' I'$) .
\end{itemize}}}
\end{remark}
\vskip 5mm We are ready to state the main result of this section.
\vskip 3mm

\begin{theorem} \label{main-transform} Let $I$ be a contracted ideal of a local regular ring $(R, \M) $
of dimension two.  With the above notation we have
$$\depth\gr_{I}(R)= {\rm min} \{ \depth\gr_{{I' }_N}(S_N) \   :   \ N \in {\rm Max}(I') \}.$$
\end{theorem}

\begin{proof} First we prove that $\gr_{I}(R)$ is Cohen-Macaulay if and only if $
\gr_{{I'}_N}(S_N)  $ is Cohen-Macaulay for every $\ N \in {\rm Max}(I').$

Assume    that  $\gr_{I}(R)$ is  Cohen-Macaulay and let $(a,b) $ be a minimal reduction of $I. $   By
Proposition \ref{dim2} we have $I ^2=(a,b)I $, and in particular, $I^2S=(a,b)IS$. By definition of the
transform of an ideal one has that
$\ell^{2d}(I')^2=\ell^{2d}(a,b)'I'=\ell^{2d}(a',b')I'.$ Since $\ell$ is regular in $S,$ it follows that
$(I')^2=(a',b')I'$ in $S,$ hence  $ (I')^2=(a',b')I'$ in $T.$ Thus $ \gr_{{I' }_N}(S_N)  $ is
Cohen-Macaulay for every $\ N \in {\rm Max}(I').$

For the converse let $(a,b) $ be  a minimal reduction of $I.$ Since
$\{a'=a/ \ell^d,\ b'=b/\ell^d\}  $ generates a minimal reduction of
$I'$ in $T, $ by Proposition \ref{dim2} we get  $(I')^2T=(a',b')I'T. $ Now, by Remark \ref{reduction}
(4), we note that ${\rm Max}(I')={\rm Max}((a',b')I').$  Hence the equality  is local on all maximal
ideals
$N$ which contain $(a',b')I', $  thus in $S$ $$(I')^2=(a',b')I'.$$ It follows that $\ell^{2d} (I')^2 =
\ell^{2d}(a',b')I' =
\ell^{d}(a',b')\ell^{d}I'$, that is, $$I^2S=(a,b)IS.$$ By Remark
\ref{reduction} (1) and (3),  both $I^2$ and $(a,b)I $ are contracted from $S, $    hence 
$I^2=I^2S\cap R=(a,b)IS\cap R=(a,b)I
$.  Therefore  $I^2=(a,b)I$  and thus $\gr_{I}(R)$ is Cohen-Macaulay. This concludes the proof of the
first part of the theorem.

Now it is enough to prove that $\depth\gr_{I}(R)>0$ if and only if
$\depth\gr_{{I' }_N}(S_N) >0 $ for every $ N \in {\rm Max}(I'). $ Assume  $\depth\gr_{I}(R)>0$. In
particular, one has that $I^{n+1}:_R a=I^n$ for every $n\ge 0$ with $a $ superficial for $I.$ Let
$a'=a/\ell^d, $   it is enough to prove that   $ (I')^{n+1}:_{S} a' =(I')^{n }$ for every $n. $    In
fact from this it follows that $a'$ is regular in $ \gr_{{I' }_N}(S_N) $ for every localization $T=S_N $
because $ (I')^{n+1}T :_T a' = ( (I')^{n+1}:_{S} a')_N = (I')^{n }_N= (I')^{n }T$ for every $n.$

\noindent  Let $c/\ell^s$ be any element of $S$, with $c\in\M^s$. Suppose   $ c/\ell^s$ is in
$(I')^{n+1}:_{S} a'$, that is,
$\frac{c}{\ell^s}\frac{a}{\ell^d}\in (I')^{n+1}$.  To prove that
$c/\ell^s\in (I')^{n}$ we distinguish two cases.

If $s\le dn$, then $$\ell^{d(n+1)}\frac{ca}{\ell^{s+d}}\in
\ell^{d(n+1)}(I')^{n+1}=I^{n+1}S,$$ that is, $\ell^{dn-s}ca\in I^{n+1}S\cap R=I^{n+1}$. Thus 
$\ell^{dn-s}c\in I^n$, and
$$\ell^{dn-s}c=\ell^{dn}\frac{c}{\ell^{s}} \in I^nS=\ell^{dn}(I')^{n}.$$  Since $\ell$ is regular in
$S$, it follows that $c/\ell^s\in (I')^n$, and this concludes this case.

Assume now that $s>dn $ and let $\M=(x,\ell).$  Since $\ell$ is coprime for $\M^{s-dn}I^{n}  $  there
exists $f \in \M^{s-dn}I^{ n} $ such that $f=x^s -p\ell $ with $p \in \M^{s-1}.$ Hence we get that
$$\left(\frac {x}{\ell}\right)^{s}= \left(\frac {p}{\ell^{s-1}}\right) + \frac{f}{\ell^s} \mbox{\ \
with \ \ } \frac{f}{\ell^s}\in (I')^n. $$ Since $c\in\M^s,$ we may write $c= ux^s+q\ell$ with $u\in R,
$ $q \in
\M^{s-1}.$ Now $$\frac{c}{\ell^s}= u\left(\frac {x}{\ell}\right)^s +
\frac{q}{\ell^{s-1}}= \frac{q+pu}{\ell^{s-1}} +  \frac{f}{\ell^s}.$$ Hence $\frac{q+pu}{\ell^{s-1}} \in
(I')^{n+1} : a' $ and we have to prove that $\frac{q+pu}{\ell^{s-1}} \in (I')^n.$ By repeating this
argument, after $s-dn $ steps we are in the already discussed case
$s<dn$.

We now assume   that $\depth\gr_{{I' }_N}(S_N) >0 $ for every $ N \in {\rm Max}(I') $ and we have to
prove that $\depth \gr_I(R) >0. $  We recall that if $J=(a,b) $ is a minimal reduction of $I,$ then
$J'=(a/\ell^d,b/\ell^d) $ is a minimal reduction of $I'$ in $T=S_N $ and, by the assumption, we may
suppose that $a'=a/\ell^d$ is   regular in $\gr_{{I' }_N}(T).$

\noindent Thus $(I')_N^{n+1}  :_T  a'=(I')_N^n  $ for every $n\ge 0  $ and for every $ N \in {\rm
Max}(I') $  which implies $(I')^{n+1} :_{S}  a'=(I') ^n  $ because it is a  local fact on the maximal
ideals
$N \in {\rm Max}(I').$

\noindent We conclude if we  prove  that $I^{n+1}:a = I^n $ for every
$n\ge 0.$

Let $b\in I^{n+1}:_R a$, that is, $ba\in I^{n+1}$. Since
$o(I^{n+1})=d(n+1)$ and $o(a)=d$, one has $o(b)\ge dn$. Then $$ ba =
\ell^{d(n+1)}\frac{ba}{ \ell^{d(n+1)}}=
\ell^{d(n+1)}\frac{b}{\ell^{dn}}\frac{a}{\ell^d}\in I^{n+1} =\ell^{d(n+1)}(I')^{n+1} ,$$ and since
$\ell$ is regular in $S $,
$\frac{b}{\ell^{dn}}\frac{a}{\ell^d}\in (I')^{n+1} $, thus
$\frac{b}{\ell^{dn}}\in (I')^n .  $ It follows that  $ b =
\ell^{dn}\frac{b}{\ell^{dn}}\in \ell^{dn}(I')^n=I^nS$, and then $b\in I^nS \cap R=I^n$, since $I^n$ is
contracted from $S.$
\end{proof}

Theorem \ref{main-transform} can be applied also in the graded setting by localizing. We present now
some corollaries which hold both in the local and in the graded case.

\begin{corollary} \label{L2} Let $I$ and $J$ be  contracted ideals with coprime characteristic forms. 
Then $$\depth\gr_{IJ}(R)= \min \{\depth\gr_{I}(R),
\depth\gr_{J}(R) \}.$$
\end{corollary}
\begin{proof} Note  that if $g=g_1^{\beta_1}\dots g_k^{\beta_k}$  is a factorization in irreducible
factors of the characteristic form $g$ of an ideal $I, $  then $$ \Max (I')=\{ (g_i/\ell,\ell)\ :\
i=1,\dots,k
\}, $$ where $\ell$ is a coprime element for $I.$ Since the characteristic forms of $I$ and $J$ are
coprime, $\Max(I') \cap
\Max(J')=\emptyset.$ Thus $(IJ)'_N=I'_N $ for every $N\in \Max(I')$ and  $(IJ)'_M=J'_M $ for every
$M\in \Max(J').$  By using twice Theorem \ref{main-transform}, we have
$$\begin{array}{rl}
  \depth\gr_{IJ}(R)& = \min \{\depth\gr_{I'_N}(S_N),
  \depth\gr_{J'_M}(S_M) \ : \ N\in \Max(I'), M \in \Max(J') \}\\
  & = \min \{\depth\gr_{I}(R), \depth\gr_{J}(R) \}.
\end{array}$$
\end{proof}

In particular:  

\begin{corollary} \label{L1} Consider the factorization   of a contracted ideal $I $ as in  Theorem
\ref{Zariski2}. Then 
$$\depth\gr_{I }(R)= \min \{\depth\gr_{L_i}(R)\ : \ i=1,\dots,k  \}.$$
\end{corollary}

\vskip 3mm

The above result leads to study the depth of the associated graded ring to a lex-segment ideal. This
will be the topic of the next  section.

\begin{remark} {\rm  Trung and Hoa  gave in \cite{TH} a combinatorial characterization of the
Cohen-Macaulayness of semigroup  rings which
 can be applied  to the study of  the Cohen-Macaulay property of the Rees algebra of monomial ideals. 
In  principle their   result  in connection with Corollary \ref{L1} can be used to give  combinatorial 
description  of the  Cohen-Macaulayness  of the associated graded rings to contracted ideals. In
practice,  however, we have not been able to obtain such a characterization.

}
\end{remark}

By Corollary \ref{L1} and Example \ref{differenze}, we have
\begin{corollary} \label{L3} Consider the factorization of a contracted ideal $I$ as in  Theorem
\ref{Zariski2}. If  $ o(L_i) \le 2 $ for every $i=1,\dots,k, $ then $ \gr_{I }(R) $ is Cohen-Macaulay.
\end{corollary}

 Consider the Rees algebra $\Rees(I)=\oplus_{n\in\N}I^n$   of $I$, and the fiber cone
$\FC(I)=\gr_I(R)\otimes R/\m$  of $I$. In the special case $o(L_i)=1,$ one has the following

\begin{theorem} \label{gcdsqfree} Let $I\subset R=k[x,y]$ be a homogeneous contracted ideal with
$o(I)=d$. Assume that the characteristic form is a square-free polynomial (it has no multiple factors).
Then: 
\begin{itemize} 
  \item[(1)] The Rees algebra  $\Rees(I)$ is a   Cohen-Macaulay normal
    domain and   the defining ideal  $J$ of   $\Rees(I)$  has  the
    expected form  in the sense of \cite[\S 8.2]{Vas} and \cite[1.2]{MU},
    that is, $J$ is the ideal of $2$-minors of a $2\times (d+1)$
    matrix $H$. 

  \item[(2)] The associated graded ring  $\gr_I(R)$  is
    Cohen-Macaulay  with Hilbert series 
    $$\HS_I(z)=\frac{\lambda(R/I)+{ d \choose 2}z}{(1-z)^2}.$$

  \item[(3)]  The fiber cone $F(I)$ is a Cohen-Macaulay reduced ring
     defined by the $2$-minors of a  $2\times d$ matrix of linear
    forms.  Furthermore $F(I)$ is a domain if and only if $I=\m^d$.
\end{itemize} 
\end{theorem} 

\begin{proof}  By Corollary \ref{L3}, $\gr_I(R) $ is Cohen-Macaulay.  Hence
$\Rees(I)$  is also  Cohen-Macaulay. Since $\GCD(I_d)$ is square free,
$I$ is given by the $d$-minors of the matrix $\phi$ of Proposition
\ref{newcontracted} with $\alpha_i\neq \alpha_j$ if $i\neq j$. It follows immediately that the ideal of
the entries of  $\phi$ is $\m$ and that the ideal of the $(d-1)$-minors of
$\phi$ is $\m^{d-1}$. By   \cite[1.2]{MU} we conclude  that $J$ has the expected form, that is, it is
given by the $2$-minors of a certain matrix $H$. We can   write down explicitly the matrix. If we
present $\Rees(I)$ as
$R[t_0,\dots,t_d]/J$  by sending $t_i$ to $(-1)^{i}$ times the
$d$-minor of $\phi$ obtained by deleting the $i+1$-th column, then $J$ is generated by the $2$-minors
of the matrix:
$$ H=\left(
\begin{array}{ccccccccccc} 
  x  &  \alpha_1t_1+y^{b_1-1}t_0  &  \alpha_2t_2+y^{b_2-1}t_1 & \dots &   
  \alpha_dt_d+y^{b_d-1}t_{d-1}  \\ 
  -y &       t_1   &t_2 &  \dots &  t_d       
\end{array}
\right). 
$$ 

By Theorem \ref{Zariski2}  $I$ is a product of    complete intersections  of order $1.$ Thus $I$ is
integrally closed. In a two-dimensional regular ring, this is equivalent to the normality of
$\Rees(I).$  This conclude the proof of part (1).   For part (3) one notes that the defining equation
of $F(I)$ are the $2$-minors of the matrix obtained from $H$ by replacing $x$ and $y$ with $0$. The
dimension of $F(I)$ is $2$ and the codimension of $F(I)$ is
$\mu(I)-2$, i.e., $d-1$. So $F(I)$ is defined by a determinantal ideal with the expected codimension,
thus it is Cohen-Macaulay. That
$F(I)$ is reduced follows by the fact that one of the initial ideals of its defining ideal is  $(t_it_j
: 1\leq i<j \leq d)$.  Finally, if $I$ is not $\m^d$, then at least one of the $b_i$, say $b_k$, is
$>1$ and then some of the generators of the defining ideal of $F(I)$ have $t_k$ as a factor. Therefore
$F(I)$ is not a domain. Now, if $I=\m^d$, then $F(I)$ is the $d$-th Veronese algebra of $R$, hence a
domain.   It remains to prove the assertion on
$\HS_I(z)$.  Since $\gr_I(R)$ is Cohen-Macaulay, its $h$-vector has length $\leq 1$. Obviously
$h_0(I)=\lambda(R/I)$. Since the $L_i$'s are complete intersections and $\beta_i =1, $ from Proposition
\ref{hfperI}  it follows   that $h_1(I)={ d\choose 2}$.  
\end{proof} 

In the above theorem it is proved that $h_1(I)={d \choose 2}$ for contracted ideal with square free
characteristic form. In general the following inequalities hold:

\begin{proposition} \label{boundh1} Let $I\subset R=k[x,y]$ be an $\m$-primary  homogeneous ideal.  If 
$I$ is monomial or contracted,  then $h_1(I)\geq {\mu(I)-1 \choose 2}$ and
$e(I)\geq \lambda(R/I)+ {\mu(I)-1 \choose 2}$.  
\end{proposition}

\begin{proof}  In general one knows that $e(I)\geq h_0(I)+h_1(I)$, see \cite[Lemma 1]{V1}. So it is
enough to prove the first inequality.  By Proposition
\ref{hfperI}, if the inequality holds for lex-segment ideals, then it holds for contracted ideals. Thus
to conclude it is enough to prove the first inequality for monomial ideals.

Let $I$ be a monomial ideal, say with associated column sequence
$a=(a_0,\dots, a_d)$ and differences sequence $b=(b_1,\dots,b_d)$. Now one has  $\mu(I) = \left | \{ i
| b_i > 0 \} \right | + 1.$ Suppose  that one of the $b_i$'s is $>1$, say $b_k>1$. Set
$c=(c_1,\dots,c_d)$ with $c_i=b_i$ if $i\neq k$ and $c_k=b_k-1$. Denote by $f$ the sequence whose
differences sequence is $c$, i.e.,
$f_0=0$ and $f_i=\sum_{j=1}^i c_j$, and by $J$ the corresponding monomial ideal. In other words,
$f_j=a_j$ if $j<k$ and $f_j=a_j-1$ if
$j\geq k$. We claim that $$h_1(I)\geq h_1(J). \eqno(2)$$ To prove this note first that
$$\lambda(R/I)-\lambda(R/J)=\sum_i a_i-\sum_i f_i=d-k+1.$$ Therefore $$h_1(I)-h_1(J)= \lambda(R/I^2) -
\lambda(R/J^2)-3 (d-k+1)$$ and hence  (2) is equivalent to:
$$\lambda(R/I^2)-\lambda(R/J^2)\geq 3 (d-k+1).$$ Denote by $a^{(2)}$ and $f^{(2)}$ the column sequences
associated with $I^2$ and $J^2$ respectively. Note that $a^{(2)}_i=a_j+a_h$ for some $j$ and $h$ with
$0\leq j,h \leq d$ and $j+h=i$.   If $i\geq 2k-1$, then at least one among $j,h$ is $\geq k$ and if 
$i\geq k+d$, then both $j,h$ are $\geq k$.  It follows that:

$$ a^{(2)}_i\geq \left\{ 
\begin{array}{ll}
  f^{(2)}_i   & \mbox{ if  } i=0,\dots, 2k-2,\\
  f^{(2)}_i+1 & \mbox{ if  } i=2k-1,\dots, k+d-1,\\
  f^{(2)}_i+2 & \mbox{ if  } i=k+d,\dots, 2d.
\end{array}
\right.
$$

We may conclude that: 
$$\lambda(R/I^2)-\lambda(R/J^2)=\sum_{i=0}^{2d} a^{(2)}_i-\sum_{i=0}^{2d} f^{(2)}_i\geq  3(d-k+1)$$ as
desired.  Since the number of generators of $I$ and $J$ is, by construction, the same, it is now enough
to prove  the assertion for
$J$.  Repeating  the argument it is enough to prove the statement for a monomial ideal $H$ whose  
differences sequence consists only of
$0$'s  and $1$'s.  Such an ideal has $\alpha+1$ generators and one generator, namely $y^\alpha$, of
degree $\alpha$. In particular, it is a lex-segment ideal with respect to $y$, whose differences
sequence does not contain 0.  After exchanging $x$ and $y$ and by applying the same procedure as above
to $H$, one ends up with a power of the maximal ideal, for which it is easy to see that the inequality
holds.
\end{proof}  

One may wonder whether the inequality $h_1(I)\geq {\mu(I)-1 \choose 2}$ holds more generally for every
$\M$-primary ideal $I$. We believe that this is indeed the case.

In  general $h_2(I) $ need not be non-negative  for an $\m$-primary ideal $I$. The ideal $I$ generated
by  $4$ generic polynomials of degree $7$ and one generic polynomial of degree $8$ (take for example
$ x^7, y^7, x^3y^4, x^6y - xy^6,  x^2y^6-x^5y^3 $) has $h_2(I) =-1$. On the other hand, there is some
computational evidence that 

\begin{conjecture} For a contracted ideal $I$ one has $h_2(I) \geq 0$.
\end{conjecture} 

Note that, in view of  Proposition \ref{hfperI}, to prove the conjecture one may assume right away that
$I$ is a lex-segment ideal. 
\bigskip

\section{Lex-segment ideals and depth of the associated graded ring}
\label{lexsegmento}

In this section we  study the depth of the associated graded ring to a lex-segment ideal in $k[x,y].$
This  is strongly motivated by Corollary \ref{L1} which moves the computation of  the depth of the
associated graded ring from   contracted ideals to lex-segment ideals. 

We start by giving  classes of lex-segment ideals whose associated graded ring has positive depth or is
Cohen-Macaulay.  Notice that  we can apply Theorem \ref{depth} to such  classes since, in our
setting,   $\gin(I)$ is a lex-segment ideal.  In the second part of the section we    find  new classes
of lex-segment ideals whose associated graded ring is Cohen-Macaulay by interpreting   Theorem
\ref{main-transform} in the case of lex-segment ideals.

\vskip 3mm Let $L$ be a lex-segment ideal in $R=k[x,y] $. As we have already seen one has $$L
=(x^d,x^{d-1}y^{a_1},x^{d-2}y^{a_2},\dots,y^{a_d}),$$ with
$0=a_0<a_1<a_2<\cdots<a_d. $ The sequence $(b_1,\dots,b_d)$,  with
$b_i=a_i-a_{i-1}$, is the differences sequence of $L$. From now on we may assume $a_d>d$, otherwise
$L=\M^d $ and its associated graded ring is Cohen-Macaulay.

By Proposition \ref{dim2}, if $I$ is an $\M$-primary ideal in $R$ with
$I^2 = JI$ for a minimal reduction $J$ of $I$, then $\gr_I(R)$ is Cohen-Macaulay. We show now that in
the class of lex-segment ideals, $L^2 = JL$ for certain kind of (non-minimal) reduction $J$, will yield
positive depth for the associated graded ring.

\begin{proposition}\label{redndepth} Let $L = (x^d, \ldots,  y^{a_d})$ be a lex-segment ideal. If $L^2 =
(x^d, x^{d-i}y^{a_i}, y^{a_d})L$ for some $i = 0, \ldots, d$, then
$\depth \gr_L(R) > 0$.
\end{proposition}

\begin{proof} We show that $L^n : (x^d, y^{a_d}) = L^{n-1}$ for all $n \geq 1$. For
$n = 1$, it is obvious. Let $n > 1$ and assume that the result is true for $n-1$. Let $h \in L^n :
(x^d, y^{a_d}).$ Without loss of generality we may assume that $h$ is a monomial. Since $L^n =
J^{n-1}L$, we may write
$$h x^d  =  (x^d)^{r_1}(x^{d-i}y^{a_i})^{r_2}(y^{a_d})^{r_3} g_1,
\eqno{(3)}$$
$$h y^{a_d}  =  (x^d)^{s_1}(x^{d-i}y^{a_i})^{s_2}(y^{a_d})^{s_3} g_2
\eqno{(4)}$$ for some $g_1, g_2 \in L$ and $\sum_i r_i = n-1 = \sum_j s_j.$ We need to show that $h \in
L^{n-1}$. If $r_1 > 0$ or $s_3 > 0$, then clearly
$h \in L^{n-1}$. Suppose $r_1 = s_3 = 0$. If $s_2 = 0$, then $s_1 = n-1$. Therefore, $x-\deg h \geq
(n-1)d$ so that $h \in L^{n-1}$. Similarly, if $r_2 = 0$, then $r_3 = n-1$. Hence $y-\deg h \geq
(n-1)a_d$ so that $h \in L^{n-1}$. Suppose $r_2 \geq 1$ and $s_2 \geq 1$. Then from (3) it follows that
$y-\deg h \geq a_{d-i}$ and from (4) it follows that $x-\deg h \geq d-i$. Therefore, $x^{d-i}y^{a_i}$
divides $h$. Write $h = x^{d-i}y^{a_i} h_1 $. Then we have
\begin{eqnarray*}
  h_1  x^d & = & 
  (x^d)^{r_1}(x^{d-i}y^{a_i})^{r_2-1}(y^{a_d})^{r_3} g_1, \\
  h_1  y^{a_d} & = & 
  (x^d)^{s_1}(x^{d-i}y^{a_i})^{s_2-1}(y^{a_d})^{s_3} g_2.
\end{eqnarray*} Therefore $h_1  \in L^{n-1} : (x^d, y^{a_d}) = L^{n-2}$, by induction hypothesis. Hence
$h = x^{d-i}y^{a_i} h_1  \in L^{n-1}$. Therefore
$L^n : (x^d, y^{a_d}) = L^{n-1}$ for all $n \geq 1$ and hence $\depth
\gr_L(R) > 0.$
\end{proof}

The following proposition gives a class of lex-segment ideals to which one can apply Proposition
\ref{redndepth}:

\begin{proposition}\label{redn1} Let $L = (x^d,  \ldots, y^{a_d})$ be a lex-segment ideal such that
$b_2 \geq b_3 \geq \cdots \geq b_d$. Then $L^2 = (x^d, x^{d-1}y^{a_1}, y^{a_d})L.$ In particular,
$\depth \gr_L(R) > 0$.
\end{proposition}

\begin{proof} Set $J = (x^d, x^{d-1}y^{a_1}, y^{a_d})$. We need to show that for all
$0 \leq i \leq j \leq d$, $x^{d-i}y^{a_i} x^{d-j}y^{a_j} \in JL.$
 
We split the proof into two cases:
\vskip 2mm
\noindent {\it Case I:} If $i+j-1 \leq d$, then we show that
$x^{2d-i-j}y^{a_i+a_j} = x^{d-1}y^{a_1}\cdot x^{d-i-j+1}y^{a_{i+j-1}}\cdot m$ for some monomial $m$.
\vskip 2mm
\noindent Consider the following equations:
\begin{enumerate}
  \item $a_{i+j-1}-a_j = b_{i+j-1}+b_{i+j-2}+\cdots+b_{j+1}$
  \item $a_i - a_1 = b_i + b_{i-1} + \cdots + b_2$.
\end{enumerate} Since $b_2 \geq b_3 \geq \cdots \geq b_d$, $a_i - a_1 \geq a_{i+j-1}-a_j$. Therefore,
$a_i + a_j \geq a_{i+j-1} + a_1$. Hence, we may write $x^{2d-i-j}y^{a_i+a_j} = x^{d-1}y^{a_1}\cdot
x^{d-i-j+1}y^{a_{i+j-1}}\cdot m$ for some monomial $m$, so that
$x^{2d-i-j}y^{a_i+a_j} \in JL$.
\vskip 2mm
\noindent {\it Case II:} If $i+j-1 > d$, then  $x^{2d-i-j}y^{a_i+a_j} = x^{d-k-1}y^{a_{k+1}}\cdot
y^{a_d}\cdot m'$ for some monomial $m'$, where $k = i+j-1-d$.
\vskip 2mm
\noindent As in {\it Case I}, write $a_d - a_i$ and $a_j - a_{k+1}$ as sum of
$b_l$'s and conclude that $a_i + a_j \geq a_d + a_{k+1}$. Therefore
$x^{2d-i-j}y^{a_i+a_j} = x^{d-k-1}y^{a_{k+1}}\cdot y^{a_d}\cdot m'$, so that $x^{2d-i-j}y^{a_i+a_j} \in
JL$.
\vskip 2mm
\noindent Therefore, for all $0 \leq i \leq j \leq d$, $x^{2d-i-j}y^{a_i+a_j}
\in JL$ and hence $L^2 = JL$. Now using Proposition \ref{redn1} we may conclude that $\depth \gr_L(R) >
0$.
\end{proof}

We apply now the theory developed in Section \ref{contracted} to
 lex-segment ideals.  Recall that a lex-segment ideal
$L=(x^d,x^{d-1}y^{a_1},\dots,y^{a_d})$ is  contracted from $S=R[\M/y]=R[x/y].$ 

In particular, $LS \cap R = L  $ and $ LS=y^dL',$ where $L' $ is the monomial ideal of $S$ generated by
the elements $  \left(
\frac{x}{y}\right)^{d-i} y^{a_i-i}  $ for every $i=0,\dots,d$.

It will be useful to consider  $\varphi :S=R[x/y]=k[x,y,z]/(x-yz)  \to P=k[y,z] $ the natural  ring
homomorphism  defined by sending the class of $f(x,y,z)$ to $f(yz,y,z)$ for every $f(x,y,z)\in
k[x,y,z].$ It is easy to see that  $\varphi$ is  an isomorphism. We set $ T(L)=\varphi(L')$, so that
$T(L)$ could  be identified with a monomial ideal in $k[x,y]$.   Moreover, one has
$$\gr_{T(L)}(P) \simeq \gr_{L'}(S).$$  These facts hold for every contracted ideal, not only for
lex-segment ideals.  In practice the ideal $T(L)$ can be obtained from $L$ by substituting $x$ with
$yz$ and dividing any generator by $y^d,$ where
$d=o(L).$ In the following examples we explain in details the procedure.

\begin{example} \label{ex1} {\rm {Let $L= (x^4,x^3y,x^2y^3,xy^4,y^{10} )\subset R=k[ x,y ]$. Note that
$o(L)=4$. The transform $L'$ of $L$ is defined to be the ideal of
$S=R[z]/(x-zy)$ such that $LS=y^4L'$. Thus one has $$LS= (y^4z^4, y^4z^3,y^5z^2,y^5z,y^{10})S =
y^4(z^4,z^3,yz^2,yz,y^{6})S=(z^3,yz,y^{6})S $$

Thus via the isomorphism $\varphi$ one has $T(L)= (z^3,yz,y^{6})$ in
$P=k[y,z].$ }}
\end{example}
\bigskip

We remark that in the particular case of a lex-segment ideal $L, $ its transform $L'$ is a primary
ideal for $N=(y,x/y)  $  or equivalently
$T(L) $ is a primary ideal for $(y,z).$ Hence, by Remark \ref{casi},
$$ \gr_{L'}(S) \simeq \gr_{L'_N}(S_N).$$

\vskip 4mm

Now we may rephrase Theorem \ref{main-transform} in the case of a lex-segment ideal.  As a consequence
one has that  to compute  the depth of the associated graded ring to $L$ one can pass to  the transform
$T(L)$ of $L$, which is in general easier to study. In particular, $\mu(T(I))\leq \mu(I)$ and $e(T(I))<
e(I)$, see
\cite[3.6]{H}. 

\vskip 5mm

\begin{theorem} \label{lex-segment} Let $L$ be a lex-segment ideal in $R=k[x,y].$ With the above
notation we have $$\depth\gr_{L}(R)= \depth\gr_{T(L)}(P).$$
\end{theorem}

\begin{proof} Since $\gr_{T(L)}(P)\simeq\gr_{L'}(S)$, it suffices to prove that
$\depth\gr_{L}(R)= \depth\gr_{L'}(S).$ The ideal $L$ is primary for
$\M=(x,y) $  and  $L' $ is primary for the maximal ideal $N=\M +(x/y),
$ hence by Remark \ref{casi}, $\gr_L(R) \simeq \gr_{L_{\M}}(R_{\M}) $ and $ \gr_{L'}(S) \simeq
\gr_{L'_{N}}(S_{N}).$ Now the result follows by using Theorem \ref{main-transform}.
\end{proof}

\medskip

As an immediate application of the theorem, we find  classes of lex-segment ideals, whose associated
graded ring is Cohen-Macaulay.

First we want  to give an explicit description of the ideal $T(L),$ where $L$ is a lex-segment ideal.  
Let $L=(B_1,\dots,B_{s+1})$ be the decomposition of the minimal set of generators of $L$ in subsets of
elements of the same degree, that is, $B_i$ is the block of the elements of degree $d+i-1$. Assume
$B_{s+1}\not=\emptyset$. In the Example \ref{ex1} one has $B_1=\{x^4,x^3y,x^2y^2,xy^3\},\
B_2=B_3=\cdots =B_6=\emptyset, \ B_{7}=\{y^{10}\}$.

\vskip 5mm Let $$|B_1|=p_1+1, \ \ \ |B_i|=p_i \ \ \  {\rm for }\ \ \ i=2,\ldots,s+1.$$
\vskip 5mm

\begin{proposition} \label{transform} With the above notation one has:
$$T(L)=(z^{d-p_1})+(z^{d-(p_1+\cdots+p_i)}y^{i-1} : i\ge 2,\  p_i\not= 0).$$
\end{proposition}
\begin{proof} By definition one has $B_1=\{x^d,\dots,x^{d-p_1}y^{p_1}\}$ and
$$B_i=\{x^{d-(p_1+\cdots+ p_{i-1}+1)}y^{p_1+\cdots +p_{i-1}+i},\dots,x^{d-(p_1+\cdots+
p_{i})}y^{p_1+\cdots+ p_{i}+i-1}\}$$ for $i\ge 2$. Applying the transform to the elements of
$B_i$ one obtains the set $T(B_i)=\{z^{d-(p_1+\cdots+ p_{i-1}+1)}y^{i-1},\dots,z^{d-(p_1+\cdots+
p_{i})}y^{i-1}\}$, hence the ideal $(T(B_i))=(z^{d-(p_1+\cdots+ p_{i})}y^{i-1})$ and this concludes the
proof.
\end{proof}

\medskip

It is natural to ask under which conditions is $T(L)$ a lex-segment ideal. As an easy consequence of 
Proposition \ref{transform} one gets a characterization:

\begin{lemma}\label{lexcharacter} Let $L$ be a lex-segment ideal. Then $T(L)$ is a lex-segment ideal if
and only if one of the following holds:
\begin{itemize}
  \item[(1)] $p_i\le 1$ for every $i\ge 2$;
  \item[(2)] $p_i\not=0$ for every $i\ge 2$.
\end{itemize} Moreover, in case $(2)$ $T(L)$ is a lex-segment ideal with respect to
$y$  and its differences sequence is $(p_{s+1},p_s, \dots,p_3,p_2)$.
\end{lemma}
\begin{proof} By Proposition \ref{transform} one has $$T(L)=(z^{d-p_1}, z^{d-(p_1+p_2)}y,
z^{d-(p_1+p_2+p_3)}y^2,\dots, z^{d-(p_1+\cdots+p_{s})}y^{s-1}, y^{s}).$$

It is clear that $T(L)$ is a lex-segment ideal with respect to $z$ if and only if $p_i\le 1$ for every
$i\ge 2$.  Note that by definition
$\sum_{i=1}^{s+1}p_i=d$ and $s+1=a_d-d+1$, and let  rewrite $T(L)$ as
$$T(L)=( y^{a_d-d}, y^{a_d-d-1}z^{p_{s+1}},\dots,y^2z^{p_4+\cdots +p_{s+1}},yz^{p_3+\cdots+p_{s+1}},
z^{p_2+\cdots+p_{s+1}}).$$

It follows that $T(L)$ is a lex-segment ideal with respect to $y$ if and only if $p_i\not=0$ for every
$i\ge 2$. When this is the case  the differences sequence is $( p_{s+1},p_s,\dots,p_3,p_2)$.
\end{proof}
\vskip 2mm
\begin{remark} {\rm {Recall that the differences sequence of $L$ is
$(b_1,\dots,b_d)$, with $b_i=a_i-a_{i-1}$. Note that the difference between the degree of
$x^{d-i}y^{a_i}$ and the degree of
$x^{d-i+1}y^{a_{i-1}}$ is $b_i-1$.  Thus the conditions of the lemma above can be written in terms of
the $b_i$'s.  In fact, condition
$(1)$ is equivalent to  $b_i\geq 2$ for $i\ge 2$, that is, the generators of $L$ of degree $>d$  have
all different degrees. Condition $(2)$ holds if and only if $b_i\in\{1,2\}$ for every $i$, that is,
there are generators in every degree between $d$ and $a_d$.}}
\end{remark}

By the above lemma, Example \ref{differenze} and Theorem
\ref{lex-segment} we get new classes of lex-segment ideals with Cohen-Macaulay associated graded ring:
\begin{proposition} \label{pi} Let $L$ be a lex-segment ideal.  Assume that one of the following holds:
\begin{itemize}
  \item[(1)]  $0<p_2\le p_3\le \dots  \le p_{s+1},$ or
  \item[(2)]  $ p_2\ge  p_3\ge \dots  \ge p_{s+1}.$
\end{itemize} Then $ \gr_{L}(R) $ is Cohen-Macaulay.
\end{proposition}

\section{Generic forms and lex-segment ideals}

Theorem \ref{depth} points to an interesting question: ``find  classes of ideals in $R$ such that the
associated graded ring of its initial ideal has positive depth".  It is known that if char $k = 0$ and
$I$ is an ideal in $R = k[x,y]$, then the generic initial ideal $\gin(I)$ is a lex-segment ideal. We
say that an ideal $I$ is a generic ideal if it is generated by generic forms of given degrees  and a
lex-segment ideal $L$ is generic, if it is the lex-segment ideal of a generic ideal. It is not always
true that the associated graded ring of a lex-segment ideal has positive depth, see Example
\ref{ex11}(a).  In this section we produce a sub-class of the lex-segment ideals, namely lex-segment
ideals of generic $\M$-primary ideals, with positive depth associated graded ring.  Let $I$ be a
generic $\M$-primary ideal in
$R$.  We begin with a lemma which will help us in identifying the structure of a generic lex-segment
ideal in $R$.  For a polynomial
$f(z) =\sum_i a_iz^i \in \ZZ[z]$, we let $\left | f(z) \right | =
\sum_i b_iz^i$ with $b_i = a_i$ if $a_0, \ldots, a_i > 0$ and $b_i = 0$ if $a_j \leq 0$ for some $j
\leq i$, and let $\Delta f(z) = \sum_i (a_i - a_{i-1})z^i$.

\begin{proposition}\label{polylemma} Let $H(z) \in \ZZ[z]$. Then
$$ H(z) = \left|\;\frac{\Pi_{i=1}^{r+1}(1-z^{d_i})}{(1-z)^2} \;\right|
$$ for some integers $d_1, \ldots, d_{r+1}$, $r \geq 1$ if and only if
$$
\Delta H(z) = 1 + z + \cdots + z^{d_1-1} - p_1z^{d_1} - p_2z^{d_1+1} - \cdots - p_sz^{d_1+s-1} - c
z^{d_1+s},
$$ where $0\leq p_1 \leq p_2 \leq \cdots \leq p_s$, $0 \leq c < p_s$ and
$\sum_{i=1}^sp_s + c = d_1$.
\end{proposition}
\begin{proof} Assume that $H(z)$ has the given form. We induct on $r$. Let $r = 1$. Then
\begin{eqnarray*}
  H(z) & = & \left |\; \frac{(1-z^{d_1})(1-z^{d_2})}{(1-z)^2} \; \right
  | \\
  & = & \left | \; (1 + z + \cdots + z^{d_1-1})(1 + z + \cdots +
  z^{d_2-1}) \right | \\
  & = & 1 + 2z + \cdots + d_1z^{d_1-1} + \cdots
  + d_1z^{d_2-1} + (d_1-1)z^{d_2} + \cdots +z^{d_1+d_2-2}.
\end{eqnarray*} Therefore $\Delta H(z) = 1+z+\cdots+z^{d_1-1}-z^{d_2}-z^{d_2+1}-\cdots z^{d_1+d_2-1}$.
Then $0 = p_1 = \cdots = p_{d_2-d_1} < 1 = p_{d_2-d_1+1} = \ldots = p_{d_2-1}.$ Also $\sum_i p_i =
(d_1+d_2-1) - (d_2-1) = d_1$. Hence the assertion follows.
\vskip 2mm
\noindent Now assume that $r > 1$ and that the assertion is true for all $l < r$. Let
$$ H'(z) = \left | \frac{\Pi_{i=1}^r(1-z^d_i)}{(1-z)^2} \right |
$$ Then by inductive hypothesis, there exist $p_1, \ldots, p_s, c$ such that $0 \leq p_1 \leq \cdots
\leq p_s; \; 0 \leq c < p_s; \; \sum_i p_i + c = d_1$ and $\Delta H'(z) = 1+z+ \cdots + z^{d_1-1} - p_1
z^{d_1} - \cdots - p_sz^{d_1+s-1} - cz^{d_1+s}$. Therefore, if we write $H'(z) = \sum_i a_iz^i$, then
$$a_i =
\left \{
\begin{array}{lll}
  i+1 \mbox{ if } i = 0, 1, \ldots, d_1-1 \\
  d_1 - \sum_{j=1}^{i-d_1+1}p_j \mbox{ if } d_1 \leq i \leq
  d_1+s-1 \\
  0  \mbox { if } i \geq d_1+s.
\end{array} \right. \eqno(5)$$ We have
$$ H(z) = \left | \Pi_{i=1}^r\frac{(1-z^{d_i})}{(1-z)^2}(1-z^{d_{r+1}})
\right |  = \left | H'(z) (1-z^{d_{r+1}}) \right |.
$$ If $d_{r+1} > \deg H'(z) = d_1+s-1$, then $\left | H'(z) (1-z^{d_{r+1}}) \right | = \left | H'(z)
\right | = H'(z)$. Therefore assume that $d_{r+1} \leq d_1+s-1$. If we set $H'(z)(1-z^{d_{r+1}}) =
\sum_i b_iz^i$, then
$$ b_i =
\left \{
\begin{array}{lll}
  a_i \mbox{ if } 0 \leq i \leq d_{r+1}-1 \\
  a_i - (j+1) \mbox{ if }  i = d_{r+1}+j.
\end{array} \right. $$ Set $h = \max\{ i \geq d_{r+1} \mid b_i > 0 \}.$ Therefore,
$$
\left | H'(z)(1-z^{d_{r+1}}) \right | = \left | \sum_i b_iz^i \right | = \sum_{i=0}^hb_iz^i =: P(z).
$$ We need to prove that $\Delta P(z)$ has the required properties. Denote by $\Delta P(z)_i$, the
coefficient of $\Delta P(z)$ in degree $i$. Then
$$
\Delta P(z)_i = \left \{
\begin{array}{ll}
\Delta H'(z)_i \mbox{ if } i \leq d_{r+1} - 1, \\
\Delta H'(z)_i-1 \mbox{ if } d_{r+1} \leq i \leq h,
\end{array}
\right.
$$ and $\Delta P(z)_{h+1} = -b_h$. Then
\begin{eqnarray*}
  \sum_{i=d_1}^h \Delta P(z)_i + b_h  & = & \sum_{i=1}^{d_{r+1}-d_1}p_i+
  \sum_{i=d_{r+1}-d_1+1}^{h-d_1+1}[p_i+1]+(a_h-h+d_{r+1}-1) \\
  & = & \sum_{i=1}^{h-d_1+1}p_i+ (h-d_{r+1}+1)+ a_h - (h-d_{r+1}+1) \\
  & = & \sum_{i=1}^{h-d_1+1}p_i + a_h\\
  & = & d_1   \; \; \;
\end{eqnarray*} because by equation $(5)$ one has  $ a_h= d_1 -
\sum_{j=1}^{h-d_1+1}p_j.$ Therefore $\Delta P(z) = \Delta H(z)$ satisfies the required properties.
\vskip 2mm
\noindent Let $\Delta H(z) = 1+z+\cdots+z^{d_1-1}-p_1z^{d_1}-\cdots-p_sz^{d_1+s-1}-cz^{d_1+s}$ with
$p_i$'s and $c$ satisfying the given properties. We prove by induction on $p_s$. Suppose $p_s = 1$.
Then $c = 0$ and $[p_1, \ldots, p_s] = [0,\ldots, 0, 1, \ldots, 1]$ for certain number of 0's, say $l$,
and 1's, say $m$. Set $d_2 = d_1+l$. Since $l+m = s$, we have
\vskip 2mm
\noindent
$$
\begin{tabular}{l|ccccccccc}
  $\deg n $ & 0 & 1 & $\ldots$ & $d_1-1$ & $d_1$ & $\ldots$ & $d_2$ &
  $\ldots$ & $d_1+d_2-1$  \\
  \hline
  $\Delta H(z)_n$ & 1 & 2 & $\ldots$ & 1 & 0 & $\ldots$ & -1 & $\ldots$
  & -1
\end{tabular} $$
\vskip 2mm
\noindent Therefore,
$$ H(z) = \left | \frac{(1-z^{d_1})(1-z^{d_2})}{(1-z)^2} \right|.
$$ Now assume that $p_s > 1$ and set $j = \max\{n \mid p_n > p_{n-1} \}$. Then we have, $1 \leq j \leq
s$ and $p_j = \cdots = p_s$. Since $p_s -1 > 0$, there exist non-negative integers $q, r$ such that
$c+s-j+1 = (p_s-1)q+r$ with $0 \leq r < p_s-1$. 
 
Define a polynomial $H'(z) \in \ZZ[z]$ such that
$$
\Delta H'(z)_i = \left\{
\begin{array}{lllll}
  1 \mbox{ if } 0 \leq i
  \leq d_1-1, \\
  -p_{i+1} \mbox{ if } d_1\leq i\le d_1+j-2, \\
  -p_j+1 \mbox{ if }
  d_1+j-1 \leq i \leq
  d_1+s-1+q, \\
  -r \mbox{ if } i =
  d_1+s+q, \\
  0 \mbox{ if } i > d_1+s+q.
  \end{array} \right. $$ Then $0 \leq p_1 \leq \ldots p_{j-1} \leq p_j-1 = \cdots = p_j-1; \; 0
\leq r < p_j-1 = p_s-1$ and $\sum_{i=1}^{j-1}p_i + (s-j+1+q)(p_j-1)+r = \sum_{i=1}^sp_i - (s-j+1) +
q(p_s-1)+r = \sum_{i=1}^sp_i+c = d_1$. Therefore, by induction there exist $d_1\leq d_2 \leq \cdots
\leq d_r$ such that
$$ H'(z) = \left | \frac{\Pi_{i=1}^r(1-z^{d_i})}{(1-z)^2} \right |.
$$ We show that $H(z) = \left | H'(z)(1-z^{d_1+j-1}) \right |$.  Note that $H'(z)_i = H(z)_i$ for $i
\leq d_1+j-2$ and $H'(z)_i = H(z)_i + i - (d_1+j-2)$. Therefore $[H'(z)(1-z^{d_1+j-1})]_i = H(z)_i$ for
$i
\leq d_1+s-1$ and $H(z)_i = 0$ for $i > d_1+s-1$. To complete the proof, we need to show that
$H'(z)_{d_1+s-1}-(p_j-1)-(s-j+2) \leq 0$. Since $H'(z)_{d_1+s-1} = H(z)_{d_1+s-1}+s-j+1 = c+s-j+1$, we
have
$H'(z)_{d_1+s-1}-(p_j-1)-(s-j+2) = c-p_s < 0.$ Therefore
$$ H(z) = \left | H'(z)(1-z^{d_1+j-1}) \right |.
$$
\end{proof}
\vskip 2mm
\noindent Using the above proposition, we describe the structure of generic lex-segment ideals in $R$.
We recall that, given a homogeneous ideal
$I\subset R$, the Hilbert series  $\HS_{R/I}(z)$ of $R/I$ is defined to be $\sum_{t\ge 0}
\HF_{R/I}(t)z^t$, where
$\HF_{R/I}(t)=\dim_k(R/I)_t$ is the Hilbert Function of $R/I$. If
$\dim R/I=0$, then $\HS_{R/I}(z)$ is a polynomial.
\begin{proposition}\label{genhilb} Let $I \subseteq R$ be an ideal generated by $r\ge 2$ generic forms
of degrees $d_1, \ldots, d_r$ respectively. Let $d=\min\{d_i\}$.  Then
\begin{itemize}
  \item[(1)] the Hilbert series of $R/I$ is such that
    $$
    \Delta \HS_{R/I}(z) =
    1+z+\cdots +z^{d-1}-p_1z^{d}-\cdots -p_sz^{d+s-1}-cz^{d+s},
    $$
    with $0\leq p_1 \leq \cdots \leq p_s$, $0 \leq c < p_s$ and $\sum_i
    p_i + c = d$.
    
  \item[(2)] $\lex(I) = (x^{d}, x^{d}y^{a_1}, \ldots, y^{a_d})$ such that
    there are $p_1+1$ elements in degree $d$ and $p_i$ elements in
    degree $d+i-1$ for $i = 2, \ldots, s$ and $c$ elements in degree
    $d+s$.
\end{itemize}
\end{proposition}
    
\begin{proof}
\begin{itemize}
\item[(1)]  The Hilbert series of $R/I$ is given by $$
    \HS_{R/I}(z) = \left |
    \frac{(1-z^{d_1})\cdots(1-z^{d_r})}{(1-z)^2} \right |,$$ for a
    simple proof of this fact, see \cite[4.3]{V}. Now the
    assertion follows directly from Lemma \ref{polylemma}.
    \vskip 2mm
    \noindent
\item[(2)] Since the Lex$(I)$ and $I$ have same Hilbert function,
    the assertion follows from (1).
\end{itemize}
\end{proof}

We set the notation for the rest of the section. Let $L = (x^d, x^{d-1}y^{a_1}, \ldots, y^{a_d})$ be a
lex-segment ideal in $R = k[x,y]$, where $a_i + 1 \leq a_{i+1}$. Recall the notation set up in Section
\ref{lexsegmento}: let $L = (B_1, \ldots, B_{s+1})$ be the block  decomposition of $L$ such that $|B_1|
= p_1+1, \; |B_i| = p_i$ for $i = 2, \ldots, s+1$. For the rest of the paper, we set $$c = p_{s+1}.$$
From Proposition \ref{genhilb}, we have $0 \leq p_1 \leq
\cdots \leq p_s$ and $0 \leq c < p_s$.

Now we proceed to prove that the associated graded rings of  generic lex-segment ideals have positive
depth. Recall that there are lex-segment ideals whose associated graded rings have depth zero, see
Example \ref{ex11}.

\begin{theorem}\label{positivedepth1} Let $L$ be a generic lex-segment ideal in $R$. Then $\depth
\gr_L(R) > 0$.
\end{theorem}

\begin{proof} We split the proof into two cases, namely $p_2 = 0$ and $p_2 > 0$.
\vskip 2mm
\noindent Let $p_2 = 0$. Note that, in this case, in degree $d$, the lex-segment ideal $L$ has only one
generator, namely $x^d$.  By Proposition
\ref{redndepth}, it is enough to prove that for some $i$, $L^2 = (x^d, x^{d-i}y^{a_i}, y^d)L$.
Following the notation set up above, for $i = 1, \ldots, s+1$, let $B_i$ denote the $i$-th block of
elements of degree $d+i-1$, of the minimal generating set of $L$. Let $x^py^q$ be the last element in
the block $B_{s}$. Then
\vskip 2mm
\noindent {\sc Claim :} $L^2 = (x^d, x^py^q, y^{a_d})L$.
\vskip 2mm To prove the claim, we need to show that, for any $1\leq i \leq j \leq d$,
$x^{2d-i-j}y^{a_i+a_j} \in JL$, where $J = (x^d, x^py^q, y^{a_d})$. As in the proof of Proposition
\ref{redn1}, we split the proof of the claim into different cases.
\vskip 2mm
\noindent We first show that if $i+j \leq d$, then $x^{2d-i-j}y^{a_i+a_j} = x^d\cdot
x^{d-i-j}y^{a_{i+j}}\cdot m$ for some monomial $m$.  It is enough to prove that $a_i + a_j \geq
a_{i+j}$.  Let $b_i = a_i - a_{i-1}$. Consider the following equations:
\begin{itemize}
  \item $a_{i+j} - a_j = b_{i+j} + b_{i+j-1} + \cdots + b_{j+1}.$
  \item $a_i = b_i + b_{i-1} + \cdots + b_1.$
\end{itemize}

\noindent Since $p_2 = 0, \; a_1 \geq 2$. Also note that since $p_2 \leq p_3
\leq \cdots \leq p_{s}$, the number of 2's appearing in $\{b_{i+j},
\ldots, b_{j+1}\}$ is at most the number of 2's appearing in $\{b_i,
\ldots, b_1\}$. Hence $a_i \geq a_{i+j} - a_j$. Therefore $a_i + a_j
\geq a_{i+j}$ and hence $x^{2d-i-j}y^{a_i+a_j} = x^d\cdot x^{d-i-j}y^{a_{i+j}}\cdot m$, for some
monomial $m$.

\vskip 2mm
\noindent Using similar arguments, we can show that
\begin{itemize}
  \item if $d < i+j \leq d+t$, where $t = d-p$, then
    $x^{2d-i-j}y^{a_i+a_j} = x^py^q \cdot
    x^{i+j-t}y^{a_{d-i-j+t}}\cdot m$ for some monomial $m$ and

  \item if $i+j > d+t$, then $x^{2d-i-j}y^{a_i+a_j} =
    x^{2d-i-j}y^{a_{i+j-d}}\cdot y^{a_d}\cdot n$ for some monomial $n$.
\end{itemize}

Therefore $x^{2d-i-j}y^{a_i+a_j} \in JL$ for all $0 \leq i \leq j \leq d$, so that $L^2 = JL$. Hence,
by Proposition \ref{redndepth}, $\depth
\gr_L(R) > 0$.

\vskip 2mm
\noindent Now let $p_2 \geq 1$. Note that in this case, there are minimal generators of $L$ in all
degrees from $d$ to $a_d$.  By Lemma
\ref{lexcharacter}, $T(L)$ is a lex-segment ideal with respect to $y$ in $P = k[z,y]$. Write $T(L) =
(y^{a_d-d},y^{a_d-d-1}z^c, y^{a_d-d-2}z^{c+p_s}, \ldots, yz^{c+p_s+\cdots+p_3}, z^{c+p_s+\cdots+p_2}).$
Then $b_1 = c$ and $b_i = p_{s-i+2}$, for $i = 2, \ldots, s$. Hence we have $b_2 \geq b_3 \geq \cdots
\geq b_s$. Therefore, $(T(L))^2 = (y^{s-1}, y^{s-2}z^c, z^{c+p_s+\cdots+p_2})T(L),$ by Proposition
\ref{redn1}. Therefore by Proposition \ref{redndepth}, $\depth \gr_{T(L)}(P) > 0$ and hence by Theorem
\ref{lex-segment} $\depth \gr_L(R) > 0$.
\end{proof}

\begin{example} {\rm{ (a) Let $I = (f_1, f_2, f_3)$ be a generic ideal such that $\deg f_1 = 5, \; \deg
f_2 = 7, \; \deg f_3 = 8$. Then, a computation as in the proof of Proposition \ref{polylemma}, will
give that $\Delta H = [1,1,1,1,1,0,0,-1,-2,-2].$ Therefore the corresponding lex-segment ideal is $L =
(x^5, x^4y^3,$ $x^3y^5, x^2y^6, xy^8, y^9)$. It can be seen that $L^2 = (x^5, y^9)L$ and hence
$\gr_L(R)$ is Cohen-Macaulay. In Theorem \ref{normalrees} we actually prove that the Rees algebra of
such ideals are normal.

(b)  Let $I = (f_1, f_2, f_3, f_4, f_5)$ be a generic ideal such that
$\deg f_1 = 10, \; \deg f_2 = 12, \; \deg f_3 = 13,\; \deg f_4=15,
\;\deg f_5 =15$. The corresponding lex-segment ideal is $L = (x^{10}, x^9y^3, x^8y^5,x^7y^6, x^6y^8,
x^5y^9, x^4y^{11}, x^3y^{12}, x^2y^{13}, xy^{14}, y^{16})$ and its Hilbert series is
$$\HS_{L}(z)=\frac{97+58 z+z^3}{(1-z)^2}.$$ By Proposition
\ref{HScorta} one has that $\depth\gr_{L}(R)=1$.}}
\end{example}
\vskip 2mm
\noindent For a generic lex-segment ideal $L$, we have seen that $|B_1|-1 \leq |B_2| \leq \cdots \leq
|B_{s}|$ and $|B_{s+1}| < |B_s|$, where $L = (B_1, \ldots, B_{s+1})$ is a block decomposition of $L$. 
Therefore, in terms of the number of generators in each degree, there can be an ``irregularity" in the
last block of elements. Since we have shown that the associated graded ring of generic lex-segment
ideals have positive depth, it is natural to ask, whether the associated graded ring is Cohen-Macaulay
when this ``irregularity" is removed. In the following theorem, we answer this question affirmatively.

\begin{theorem}\label{normal} Let $L$ be a generic lex-segment ideal in $R$ such that $c = 0$.  Let
$J = (x^{d-p_1}y^{a_{p_1}}, x^d+y^{a_d})$. Then $L^2 = JL$ and $L$ is integrally closed. In particular,
$\gr_L(R)$ is Cohen-Macaulay.
\end{theorem}
\begin{proof} Let $L = (x^d, x^{d-1}y^{a_1}, \ldots, xy^{a_{d-1}}, y^{a_d})$. Following our previous
notation, let $p_1+1 = | B_1 |$ and $p_i = | B_i |$ for $i = 2, \ldots, s$. First we prove that $L^2 =
JL$ for $J = (x^{d-p_1}y^{a_{p_1}}, x^d+y^{a_d})$. We show that for $1\leq i \leq j
\leq d$, $x^{2d-i-j}y^{a_i+a_j} \in JL$. We split the proof into three cases: 
\vskip 2mm
\noindent {\it Case I:} $\;\;i+j < p_1$.
\vskip 1mm
\noindent Write $x^{2d-i-j}y^{a_i+a_j} = (x^d+y^{a_d})(x^{d-i-j}y^{a_i+a_j}) -
x^{d-i-j}y^{a_i+a_j+a_d}.$ Note that for $i \leq p_1, \; a_{i-1}+1 = a_i$. Therefore, $a_i+a_j =
a_{i+j}$, if $i+j < p_1$. Hence
$(x^d+y^{a_d})x^{d-i-j}y^{a_i+a_j} \in JL$. Now,
$x^{d-i-j}y^{a_i+a_j+a_d} = (x^{d-p_1}y^{a_{p_1}}) (x^{d-(d-p_1+i+j)}y^{a_i+a_j+a_d-a_{p_1}}).$ Since
the number of minimal generators of $L$ in each degree in increasing, as argued in the proof of Theorem
\ref{positivedepth1}, we can show that
$a_i+a_j+a_d \geq a_{d-p_1+i+j}+a_{p_1}.$ Therefore
$(x^{d-(d-p_1+i+j)}y^{a_i+a_j+a_d-a_{p_1}})$ $ \in L$ so that
$x^{2d-i-j}y^{a_i+a_j} \in JL$.
\vskip 2mm
\noindent {\it Case II:} $\;\;p_1 \leq i+j \leq d+p_1$.
\vskip 1mm
\noindent Writing $a_j - a_{p_1}$ and $a_{i+j-p_1} - a_i$ as in the proof of Theorem
\ref{positivedepth1}, one can easily see that, in this case
$a_i+a_j \geq a_{p_1}+a_{i+j-p_1}$. Hence $x^{2d-i-j}y^{a_i+a_j} \in JL$.
\vskip 2mm
\noindent {\it Case III:} $\;\; d+p_1 < i+j$.
\vskip 1mm
\noindent Then $2d-i-j = d-p_1-k$ for some $k \geq 1$. Therefore we can write
$x^{2d-i-j}y^{a_i+a_j} = (x^d+y^{a_d})(x^{d-p_1-k}y^{a_i+a_j-a_d}) - x^{2d-p_1-k}y^{a_i+a_j-a_d}$.
Arguments similar to that of in the proof of {\it Case I} will show that $x^{d-p_1-k}y^{a_i+a_j-a_d} \in
L$ and $x^{2d-p_1-k}y^{a_i+a_j-a_d} \in JL$. Hence $x^{2d-i-j}y^{a_i+a_j} \in JL$. Therefore
$L^2 = JL$. 

Now we proceed to prove that $L$ is integrally closed. From  Corollary \ref{L1}, it follows that if $L$
has
$r$ generators in the initial degree, then $L = \m^rN$ for a lex-segment ideal $N$. It can easily be
seen that if $L$ is generic, then so is $N$. Note also that there is only one generator in the initial
degree (i.e., $p_1 = 0$) and $c = 0$ for $N$. We have considered such ideals in the next section. In
Theorem
\ref{normalrees} we have proved that lex-segment ideals with $p_1 = c = 0$ are integrally closed.
Therefore $N$ is integrally closed. Since
$L$ is a product of power of the maximal ideal (which is integrally closed) and $N$, $L$ is integrally
closed. Hence $\gr_L(R)$ is Cohen-Macaulay.
\end{proof}

We end the section with another interesting class of lex-segment ideals whose associated graded rings
are Cohen-Macaulay, namely lex-segment ideals corresponding to ideals generated by generic forms of
equal degree.

\begin{proposition} Let $I$ be an $\M$-primary ideal generated by generic forms of equal degree.  Let
$L$ be the lex-segment ideal corresponding to $I$. Then
$\gr_L(R)$ is Cohen-Macaulay.
\end{proposition}
\begin{proof} Let $I$ be generated by $r$ forms of degree $d$. Then the Hilbert series of $R/I$ is $$
\HS_{R/I}(z) = \left | \frac{(1-z^d)^r}{(1-z)^2}
\right |.$$ A direct computation shows that $HF_{R/I}(n) = n+1 \;\;
\mbox{ for } n = 0, \ldots, d-1$, and  for $i\ge 0$
$$ HF_{R/I}(d+i) = \left\{ 
\begin{array}{ll}
  d-(r-1)(i+1)& \mbox{ if } d-(r-1)(i+1)  \geq 0\\
  0 & \mbox{ otherwise}
\end{array}.
\right.
$$ Therefore
$$
\Delta H = [1, 1, \ldots, 1, -r+1, -r+1, \ldots, -r+1, -c],
$$ where $1$ is repeated $d-1$ times, $-r+1$ is repeated $[d/(r-1)]$ times, where $[d/(r-1)]$ denotes
the largest integer smaller or equal to $d/(r-1)$, and $0 \leq c < r-1$.  Hence the corresponding
lex-segment ideal $L$ have $r$ generators in degree $d$, $r-1$ generators in degree $d+j$ for $j = 1,
\ldots, d+[d/(r-1)]-1$ and $c$ generators in degree $d+[d/(r-1)]$.  Thus, by Proposition \ref{pi},
$\gr_L(R)$ is Cohen-Macaulay.
\end{proof}

\section{Rees algebras of lex-segment ideals}

In this section we study the Rees algebras of lex-segment ideals. For an ideal $I$ in a ring $R$, the
Rees algebra $\Rees(I)$ is defined to be the $R$-graded algebra $\oplus_{n\ge 0} I^n$. It can be
identified with the $R$-subalgebra, $R[It]$ of $R[t]$ generated by $It$, where
$t$ is an indeterminate over $R$.

Let $I=(x^d,x^{d-1}y^{a_1},\ldots,y^{a_d}) $ be a lex-segment ideal in
$R = k[x,y].$ Consider the epimorphism of $R$-graded algebras $$\psi: R[T_0,\ldots,T_d] \longrightarrow
\Rees(I)$$ defined by setting
$\psi(T_i) = x^{d-i}y^{a_i}t$ and let $H = \ker\psi$ be the ideal of the presentation of $\Rees(I)$.
The goal of this section is to describe explicitly a Gr\"obner basis of $H$, for some of the
interesting classes of lex-segment ideals we have considered. We begin by describing a set of binomials
which are not in $\ker \psi$.

\begin{lemma}
\label{relazioni} Let $ a,b,c,f,g$ be integers bigger than or equal to $0$. The ideal
$H$ does not contain non-zero elements of the following forms:
$$T_i^aT_{i+1}^b-y^cT_j^fT_{j+1}^g\mbox{ with } 1\le i,j\le d-1, $$
$$y^aT_i^bT_{i+1}^c-y^fT_j^g \mbox{ with } 0\le i,j\le d-1, $$
$$T_0^aT_d^bT_j^l-y^cT_0^fT_d^gT_{k}^h\mbox{ with } 1\le j\not=k\le d-1, 0\le l,h\le 1. $$
\end{lemma}

\begin{proof} Suppose that $m_i-m_j=T_i^aT_{i+1}^b-y^cT_j^fT_{j+1}^g$ is in $H$; then
$\psi(T_i^aT_{i+1}^b)=\psi(y^cT_j^fT_{j+1}^g)$, and by comparing the degrees respectively of $t,x$ one
has:
$$\hspace*{0.1in}\left\{
\begin{array}{l} a+b=f+g\\ a(d-i)+b(d-i-1)=f(d-j)+g(d-j-1) \\
\end{array}
\right.
\eqno(6)
$$

\vskip 2mm Since $a+b = f+g$, $ad+bd = fd+gd$ and thus from $(6)$, it follows that $ai+ bi + b = fj +
gj + g$. Therefore $g-b = (f+g)(i-j)$.

If $i > j$, then, since $f+g > 0$, $g - b > 0$. Again from $(6)$, we get that $g - b \geq f+g$, i.e.,
$-b \geq f$. Therefore, the only possibility is that $f = b = 0$ and hence $i = j+1$. This implies that
$a = g$ and $c = 0$. Hence $m_i = m_j$.

If $i<j$, then one concludes in the same way as before, since one has
$b-g=(a+b)(j-i)$.

If $i=j$, then by $(6)$ one has $b=g$, and therefore $a=f$ and $c=0$. Again we have $m_i=m_j$.

\medskip

Identical arguments will show that a non-zero equation of the form
$y^aT_i^bT_{i+1}^c-y^fT_j^g$, with $ 0\le i,j\le d-1, $ is not  in $H$.

\medskip

Suppose now that an element of the form
$m_j-m_k=T_0^aT_d^bT_j^l-y^cT_0^fT_d^gT_{k}^h$ is in $H$. By a degree comparison this implies :

$$\left\{
\begin{array}{l}
  a+b+l=f+g+h\\
  ad+l(d-j)=df+h(d-k)\\
  ba_d+la_j=c+ga_d+ha_k\
\end{array}.
\right.
$$ We distinguish different cases.

If $l=h=0$, then one has $a=f$, $b=g$, and $c=0$. Thus $m_j=m_k$.

If $l=0$ and $h=1$, then it follows that $ad=fd+d-k$, that is,
$(a-f)d=d-k$.  This is a contradiction since $d$ cannot divide $d-k$. If $l=1$ and $h=0$, one concludes
as in the previous case that
$m_j=m_k$.

If $l=h=1$, then $ad+d-j=fd+d-k$, that is, $(a-f)d=j-k$. This implies that $d$ divides $j-k$ and this
is a contradiction, since $1\le j\neq k\le d-1$.
\end{proof}

In the following two propositions, we describe explicitly a Gr\"obner basis for the presentation ideal
of Rees algebras of lex-segment ideals with increasing and decreasing differences sequence, already
considered in Example \ref{differenze}.

\begin{proposition}\label{gb1} Let $I$ be a lex-segment ideal in $R$. Suppose that its differences
sequence is such that $b_i\le b_{i+1}$ for $i=1,\dots,d-1$. Then the set of elements
$$\left\{
\begin{array}{l}
  xT_i-y^{b_i}T_{i-1}, \ \ i=1,\dots,d \\
  T_iT_{j-1}-y^{b_i-b_j}T_{i-1}T_j, i,j\in\{1,\dots,d\},\ d\ge i>j\ge1
\end{array}
\right\}
$$ form a Gr\"obner basis of $H$ with respect to any term order such that the initial term of any of
the elements above is the term on the left side. Also, the Rees algebra $\R(L)$ is normal. 
\end{proposition}
\begin{proof} Let $Q$ be the ideal generated by $xT_i,\ i=1,\dots,d$, and
$T_iT_{j-1},\  i>j\ge 1$. Since the binomial relations form a universal Gr\"obner basis of $H$, to
prove that $Q=\ini(H)$ it suffices to prove that there are no relations in $H$ involving only monomials
which are not in $Q$. Note that such monomials are of the form $x^aT_0^b,\; y^aT_i^bT_{i+1}^c, \;
y^aT_i^b$ for some $a,b,c$. Since $\psi(x^aT_0^b)$ involves only $x$, it cannot be the term of an
element in $H$.  Moreover, if $y^aT_i^b-y^cT_j^e$ is in $H$, then it is easy to see that $b=e$ and
$i=j$.  From Lemma \ref{relazioni} it follows that the given relations are a  Gr\"obner basis of $H$. 
Note that the initial terms of the elements of the Gr\"obner basis are square-free, thus by
\cite[13.15]{St} the Rees algebra  $\R(L)$ is normal. 
\end{proof}

\begin{proposition} Let $L$ be a lex-segment ideal in $R$. Suppose that its differences sequence is
such that $b_i\ge b_{i+1}$ for
$i=1,\dots,d-1$. Then the set of elements
$$\left\{
\begin{array}{l}
  xT_i-y^{b_i}T_{i-1}, \ \ i=1,\dots,d \\
  T_iT_{j-1}-y^{b_i-b_j}T_{i-1}T_j, i,j\in\{1,\dots,d\},\ 1\le i<j
  \le d,
\end{array}
\right\}
$$ form a Gr\"obner basis of $H$ with respect to any term order such that the initial term of any of
the elements above is the term on the left.
\end{proposition}

\begin{proof} Let $Q$ be the ideal generated by $xT_i,\ i=1,\dots,d$, and
$T_iT_{j-1},\ 1\le i<j\le d$.  In particular, $T_i^2\in Q$ for
$i=1,\dots,d-1$.  Thus the monomials which are not in $Q$ are the ones of the form $x^aT_0^b,\
y^aT_0^bT_d^cT_j^e$, with $0<j<d$, $0\le e\le 1$, and some $a,b,c$.  Using Proposition \ref{relazioni}
and arguing as in Proposition \ref{gb1}, one concludes that $Q=\ini(H)$.
\end{proof} It is easy to see  that in general a lex-segment ideal $L$ as in the proposition above is
not integrally closed, thus $\R(L)$ is not normal.

\vskip 2mm
\noindent  In the following theorem, we obtain the Gr\"obner basis for the presentation ideal of the
Rees algebra of another sub-class of generic lex-segment ideal and then use it to produce another class
of normal Rees algebras.

\begin{theorem}\label{normalrees} Let $L$ be a generic lex-segment ideal in $R$ such that $c = p_1 = 0$.
Then the set of elements
$$\left\{
\begin{array}{lll}
  xT_i - y^{b_i}T_{i-1}, \; \; i = 1, \ldots, d \\
  T_iT_j - y^{\alpha}T_0T_{i+j}, \; \; 1\leq i \leq j < d, \;  1
  \leq i+j \leq d \mbox{ and } \alpha = a_i+a_j-a_{i+j} \\
  T_iT_j - y^{\beta}T_{i+j-d}T_d, \;\; 1 \leq i \leq j < d, \; d
  < i+j  \mbox{ and } \beta = a_i+a_j-(a_{i+j-d}+a_d).
\end{array}
\right\}
$$ form a Gr\"obner basis for $H$ with respect to any term order such that the initial term of any of
the elements in the above set is the term on the left. Also, the Rees algebra $\R(L)$ is normal.
\end{theorem}

\begin{proof} Since $c = p_1 = 0$, it follows from Theorem \ref{normal} that $L^2 = JL$ for $J = (x^d,
y^{a_d})$. Let $\B$ denote the set of elements given in the statement of the theorem. We first show
that $\B$ is a Gr\"obner basis for $H$.
\vskip 2mm
\noindent Let $Q$ be the ideal generated by $xT_i, \; i = 1, \ldots, d$ and
$T_iT_j, \; 1 \leq i \leq j < d$. The monomials which are not in $Q$ are of the form $x^ay^bT_0^c, \;
y^aT_0^bT_d^cT_j^e, \; y^aT_i^b$ with
$0 < j < d, \; 0 \leq i \leq d$ and some non-negative integers $a, b, c, e$. Then, as in the proof of
Proposition \ref{gb1}, it can easily be seen that $\ini(H) = Q$ and hence $\B$ is a Gr\"obner basis for
$H$.
\vskip 2mm
\noindent Now we prove that $\R(I)$ is normal. Since $\R(I)$ is a semi-group ring, it is enough to
prove that for any three monomials $f, g, h \in
\R(I)$, if for some integer $p$, $f^p = g^ph$, then $h = h_1^p$ for some monomial $h_1 \in \R(I)$ (see
\cite[6.1.4]{BH}). Let $S = R[T_0,
\ldots, T_d]$. Then $\R(I) \cong S/H$. Since $\B$ is a Gr\"obner basis for $H$ and $\ini(H) = Q$, the
set
$$ Q' = \left\{
\begin{array}{ll}
  x^ay^bT_0^c, \;\; a, b, c \geq 0 \\
  y^aT_0^bT_d^cT_j^e, \;\; 0 < j < d;\; a,b,c,e \geq 0 \\
  y^aT_i^b, \; a, b\geq 0;\  0 \leq i  \leq d
\end{array}
\right\}
$$ form a monomial basis for $S/H$. Note that any monomial in $S/H$ will be a power of either of the
above forms. Let $f, g, h \in S/H$ be monomials such that
$$ f^p = g^ph \mbox{ for some } p \geq 0.
\eqno{(7)}$$

\vskip 2mm

Let $f = x^ay^bT_0^c$ for some integers $a, b, c$. Then, from (7) and comparing the $x$ and $y$ degrees
of $\psi(f), \psi(g)$ and $\psi(h)$, we can conclude that both $g$ and $h$ can not contain $T_j$ for $j
\neq 0$.  Write $g = x^{a_1}y^{b_1}T_0^{c_1}$.  From $(7),$ it follows that $a_1 \leq a, \; b_1 \leq b$
and $c_1 \leq c$.  Therefore, $h = (x^{a-a_1}y^{b-b_1}T_0^{c-c_1})^p$.  Hence $h = h_1^p$, for $h_1 =
x^{a-a_1}y^{b-b_1}T_0^{c-c_1}$.
\vskip 2mm

Now assume that $f = y^aT_0^bT_d^cT_j^e$ for some $0 < j < d$ and non-negative integers $a, b, c, e$.
If $g = x^{a_1}y^{b_1}T_0^{c_1}$. Then comparing the $x$-degrees, we get that $a_1 = 0$. Therefore $g =
y^{b_1}T_0^{c_1}$, such that $b_1 \leq a$ and $c_1 \leq b$. Therefore
$h = (y^{a-b_1}T_0^{b-c_1}T_d^{c}T_j^{e})^p$. Hence $h = h_1^p$ for
$h_1 = y^{a-b_1}T_0^{b-c_1}T_d^{c}T_j^{e}$. Suppose $g = y^{a_1}T_0^{b_1}T_d^{c_1}T_i^{e_1}$ for some
$i$. Then again from (7), it follows that $i = j$ and $h = h_1^p$ for $h_1 =
y^{a-a_1}T_0^{b_1}T_d^{c_1}T_j^{e-e_1}$.
\vskip 2mm

Let $f = y^aT_i^b$ for some $a, b \geq 0$ and $0 \leq i \leq d$.  Then it is obvious from (7)  that $g$
and $h$ have to be of the same form. Thus, as in the previous cases, we conclude that $h = h_1^p$ for
some
$h_1$.

Therefore $S/H$ is normal and hence $\R(L)$ is normal.
\end{proof}

\end{document}